\documentclass[12pt]{article}
\usepackage{graphicx,psfrag,epsfig, color,float}
\usepackage{amssymb,amsmath,amscd,amsthm}
\usepackage{graphicx,psfrag,epsfig}

\usepackage{graphicx}
\usepackage[active]{srcltx}

\newtheorem{theorem}{Theorem}[section]
\newtheorem{corollary}[theorem]{Corollary}

\newtheorem{lemma}[theorem]{Lemma}

\date{}

\begin{document}

\date{}
\title{Diffusion in the presence of cells with semi-permeable membranes}
\author{
 M. Freidlin\footnote{Dept of Mathematics, University of Maryland,
College Park, MD 20742, mif@math.umd.edu}, L.
Koralov\footnote{Dept of Mathematics, University of Maryland,
College Park, MD 20742, koralov@math.umd.edu}
%%\footnote{Dept of Mathematics, University of Maryland, College
%%Park, MD 20742, koralov@math.umd.edu}
%%\footnote{Dept of Mathematics, University of Maryland,
%%College Park, MD 20742, mif@math.umd.edu}
} \maketitle

\begin{abstract}
We consider processes that coincide with  a given diffusion process except on the boundaries of a finite collection of domains. 
The behavior on each of the boundaries is asymmetric: 
the process is much more likely to enter the interior of the domain than to enter the interior of its complement,
with a small parameter controlling the trapping mechanism.  We describe the limiting behavior of the processes. In particular, if 
the parameters controlling the boundary behavior have different orders of magnitude for different domains or if the domains are nested, metastable distributions between the trapping regions are described.
\end{abstract}

{2010 Mathematics Subject Classification Numbers: \ 60F10, 35B40, 35J25, 47D07, 60J60.}

{ Keywords: Metastability, Non-standard Boundary Problem, Asymptotic Problems for PDEs. }

\section{Introduction} \label{intro}
Consider particles diffusing in a $d$-dimensional space (we will assume that the space is just the torus
$\mathbb{T}^d$ in order to avoid the issues of recurrence-vs-transience that are not of interest in the current paper). The process governing 
the motion of a particle starting at $x$ depends on a small parameter $\varepsilon$ and is denoted by $X^{x, \varepsilon}_t$.

Let $\mathcal{D} = \{D_1,...,D_n\}$, $D_k \subset \mathbb{T}^d$, be a collection of open simply connected domains with sufficiently smooth 
boundaries $\partial D_k$, $k =1,...,n$.
The boundaries are assumed to be disjoint; they model semi-penetrable membranes for the process $X^{x, \varepsilon}_t$. In $\mathbb{T}^d \setminus \bigcup_{k=1}^n \partial D_k$,
the process coincides with a given $\varepsilon$-independent diffusion, say, a Wiener process, 
while its behavior on the membranes $\partial D_k$ is asymmetric: starting at $x \in \partial D_k$, the process ``goes to the interior of $D_k$" with probability
$1/(1 + \varepsilon_k)$ and ``goes to the exterior of $D_k$" with probability $\varepsilon_k/(1 + \varepsilon_k)$, where $0 < \varepsilon_k = \varepsilon_k(\varepsilon) \ll 1$.
Actually, since one can't define the direction of the first exit of a Wiener process from $\partial D_k$, defining $X^{x, \varepsilon}_t$ rigorously involves
specifying the generator of the process, in particular, the domain of the generator (this is done in Section~\ref{dttpA}). Alternatively, one could give
rigorous meaning to the statement that the process goes to the interior or the exterior with prescribed probabilities by approximating $X^{x, \varepsilon}_t$
with processes that experience an instantaneous jump of size $\delta$ in the direction orthogonal to $\partial D_k$  upon reaching $\partial D_k$. The jump is directed to the
interior of $D_k$ with probability $1/(1 + \varepsilon_k)$ and to the exterior of $D_k$ with probability $\varepsilon_k/(1 + \varepsilon_k)$. Upon taking $\delta \downarrow 0$, one can obtain the desired process $X^{x, \varepsilon}_t$ in the limit. 

\begin{figure}[htbp]
    \centerline{\includegraphics[height=3in, width= 3.5in,angle=0]{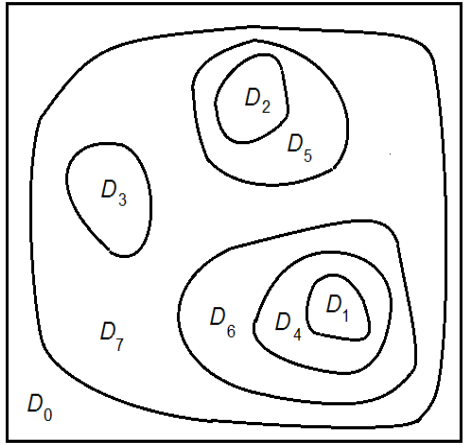}}
\vskip -15pt
    \label{secondkind}
		\caption{An example with mutiple nested domains}
\end{figure}

Our goal is to describe the behavior of $X^{x, \varepsilon}_t$ on long time intervals that grow together with $\varepsilon^{-1}$. Assume that for every pair of domains
$D_k$ and $D_l$ with $k \neq l$, either the domains are disjoint or one is a subset of another. We'll also adjoin $D_0 = \mathbb{T}^d$ to 
the list, and assume that other domains are proper subsets of $D_0$. These assumptions allow us to define the notion of the rank of $D_k$ 
inductively. We'll say that $D_k$ has rank one if it does not contain other domains. Having defined all the domains of ranks $1,..,r-1$, we'll say that
$D_k$ has rank $r$ if it is not a domain of rank that is less than $r$ and all the domains it contains have rank less than $r$. We'll write that
$D_l \prec D_k$ if ${\rm rank}(D_l) + 1 = {\rm rank}(D_k)$ and $D_l \subset D_k$. In the example shown in  Figure~\ref{secondkind}, there are three 
domains of rank one: $D_1, D_2, D_3$, two domains of rank two: $D_4, D_5$, one domain of rank three: $D_6$, and one domain of rank four: $D_7$, and
one domain or rank five: $D_0 = \mathbb{T}^d$.

\vskip -15pt

%\bb The corresponding graph $\Gamma$ is defined to have the set of vertices $\{0,...,n\}$. There is a directed edge from $l$ to $k$ if $D_l \prec D_k$.  
%The graph corresponding the the domains from \br Figure ... \be is shown in \br Figure ... \be. 
%\be

The limiting behavior of $X^{x, \varepsilon}_t$ will be described using processes with instantaneous re-distribution and
reflection, formally defined in Section~\ref{dttp}. Here, we give an intuitive description of such processes and then illustrate the asymptotic behavior of $X^{x, \varepsilon}_t$ using the example shown in Figure~\ref{secondkind}. Let $0 \leq k \leq n$ and $S \subseteq \{1,...,n\}$ be such that $D_l \prec D_k$ for each $l \in S$. 
Let $U = D_k \setminus \bigcup_{l \in S} \overline{D}_l$.
For each such $k$ and $S$, consider the corresponding process $Y^{x}_t$, $Y^{x}_0 = x \in {\overline{U}}$,  that coincides with the Wiener process in $U$, is reflected at $\partial D_k$, and  is instantaneously re-distributed on $\partial D_l$ (according to the volume measure)
upon reaching $\partial D_l$  and then 
reflected  into $U$. We stress that $Y^{x}_t$ depends on $k$ and $S$, although this is not reflected in the notation. The rigorous definition of such a process involves specifying its generator (in Section~\ref{dttp} this is done after modifying the
state space of the process in order to make the trajectories continuous). Its transition probabilities 
$\mathrm{P}(Y^{x}_t \in B) = p(t,x,B)$ can be shown to
satisfy
\[
\frac{\partial p(t,x,B)}{\partial t} = \frac{1}{2} \Delta p(t,x,B),~~~p(0,x,B) = \chi_B(x),~~~x \in U,
\]
\[
(\nabla p(t,x,B), n(x)) = 0,~~~x \in \partial D_k,
\]
\[
\int_{\partial D_l} (\nabla p(t,x,B), n(x)) \nu_l (d x) = 0,~~~l \in S,
\]
\[
p(t,x,B) = c_l(t)~~~{\rm for}~{\rm some}~c_l(t)~{\rm and}~{\rm all}~x \in \partial D_l.
\]
Here, $n(x)$ is the unit exterior normal at $ x \in \partial U$, $\nu_l$ is the volume measure on $\partial D_l$, and the
values of $c_l(t)$ are not prescribed. Thus, solving this equation, i.e., finding $p(t,x,B)$ as a function $(t,x)$ with $B$ fixed, involves
finding its boundary values $c_l(t)$. The existence and uniqueness of solutions to such non-standard problems (or, rather, the corresponding 
elliptic problems) has been discussed, e.g., in \cite{Nov} (see also \cite{FrK1}, \cite{FrK2}, where $\nu_l$ is allowed to be an arbitrary measure, and
where the corresponding process was constructed).

% In order to make the trajectories of  $Y^{k, S, x}_t$ continuous, one can modify the state space 
%of the process.  Namely, let $U'$ be the metric space obtained from $\overline{U}$ by identifying all points of $\partial D_l$, turning every $\partial D_l$ into one %point $d_l$. Let $h: \mathbb{T}^d \rightarrow U'$ be the mapping defined by $h(x) = x$ for \br $x \in U$ and $h(x) = d_k$ for $x \in \overline{D}_k$. 
%For $x \in \overline{D}_k$, we can view $Y^{k, S, h(x)}_t$ the process with values in $U'$ that is equal to $h(Y^{k, S, h(x)}_t)$.  \be

% \bb One can allow the index $x$ (formerly, the initial point for the process $Y^{k, S, x}_t$) to belong to $D_l \prec D_k$. This is accomplished by modifying the %state space of the process

%\bb One can allow the initial point $x$ for the process $Y^{x}_t$ to belong to $\mathbb{T}^d \setminus 
%\overline{D}_k$ -- in this case the process coinsides with the Wiener process until the first time it hits $\partial D_k$, and then the process continues as
%described above. \be

Having described the processes $Y^{x}_t$, let us now use the example shown in  Figure~\ref{secondkind} to discuss the asymptotic behavior of $X^{x, \varepsilon}_{t(\varepsilon)}$. Assume that $\varepsilon_k(\varepsilon) =
\varepsilon$ for each $1 \leq k \leq 7$. The distribution of $X^{x, \varepsilon}_{t(\varepsilon)}$ depends on the initial point $x$ and the time scale $t = t(\varepsilon)$.
First, consider the case when $1 \ll t(\varepsilon) \ll \varepsilon^{-1}$. In this case, as we'll see, the process has enough time to enter one the domains of rank one, but not
enough time to exit a small neighborhood of any domain $D_k$, $k \geq 1$. Thus, if $x \in \overline{D}_k$ with $k = 1,2,3$, then the distribution of $X^{x, \varepsilon}_{t(\varepsilon)}$ will be asymptotically close to the
limiting distribution of the Wiener process with reflection on the boundary of $D_k$ (i.e., the uniform distribution on $D_k$, denoted by $\lambda_k$). If $x \in 
\overline{D}_6 \setminus \overline{D}_1$, then
the process enters  $D_1$ and, asymptotically, is distributed uniformly on $D_1$. Similarly, for $x \in \overline{D}_5 \setminus \overline{D}_2$, 
$X^{x, \varepsilon}_{t(\varepsilon)}$ tends to 
a uniform distribution on $D_2$. 

The situation is slightly more complicated if $x \in \overline{D}_7 \setminus (\overline{D}_6 \bigcup \overline{D}_5 \bigcup \overline{D}_3)$. 
In this case, we need to consider the Wiener process with 
reflection on $\partial D_7$ ($Y^{x}_t$ corresponding to $k = 7$ and $S = \O $ in the above notation). Let
$\tau^{x, \varepsilon} = \tau^{x, \varepsilon}(\partial D_6 \bigcup \partial D_5 \bigcup \partial D_3)$ be the first time when the process  $X^{x, \varepsilon}_{t}$
hits $\partial D_6 \bigcup \partial D_5 \bigcup \partial D_3$. We define $\tau^x$ similarly, but for the process $Y^x_t$ rather than $X^{x,\varepsilon}_t$. If $X^{x,\varepsilon}_t$ reaches $\partial D_6$ first, then it will tend to the uniform distribution on $D_1$, 
if it reaches $\partial D_5$ first, then it will tend to the uniform distribution on $D_2$, and if it reaches $\partial D_3$ first, then it will tend to the uniform
distribution on $D_3$. Using the fact that $Y^x_t$ serves as a good approximation for $X^{x,\varepsilon}_t$ at these time scales (Section~\ref{seconv2}), we will be able to conclude that
 the distribution of $X^{x, \varepsilon}_{t(\varepsilon)}$ tends to 
\begin{equation} \label{jijil}
\mathrm{P}(Y^x_{\tau^{x}} \in \partial D_6) \lambda_1 + \mathrm{P}(Y^x_{\tau^{x}} \in \partial D_5) \lambda_2 + \mathrm{P}(Y^x_{\tau^{x}} \in \partial D_3) \lambda_3.
\end{equation}
It should be pointed out that the coefficients in (\ref{jijil}) do not depend on $\varepsilon$ and can be calculated as solutions of the appropriate boundary value problems.

Next, we discuss what happens if $x \in D_0 \setminus \overline{D}_7$. In this case, we consider the Wiener process in $D_0$ ($Y^x_t$ corresponding to $k = 0$ and $S = \O$) until the first time it hits $\partial D_7$. 
Let $\pi(x, y)$, $x \in D_0 \setminus \overline{D}_7$, $y \in \partial D_7$, be the corresponding Poisson kernel. We apply the arguments above, but starting with the
point $y$ where the process first reaches $\partial D_7$. Thus the distribution of $X^{x, \varepsilon}_{t(\varepsilon)}$ tends to 
\[
\int_{\partial D_7} (\mathrm{P}(Y^y_{\tau^{y}} \in \partial D_6) \lambda_1 + \mathrm{P}(Y^y_{\tau^{y}} \in \partial D_5) \lambda_2 + 
\mathrm{P}(Y^y_{\tau^{y}} \in \partial D_3) 
\lambda_3) \pi(x, y) dy.
\]

Now let us consider the case when $\varepsilon^{-1} \ll t(\varepsilon) \ll \varepsilon^{-2}$. In this case,  $X^{x, \varepsilon}_{t}$ has enough time to exit a  domain of rank one, but not a small neighborhood of a domain of rank two. Moreover, while the process may enter and exit a domain of rank one (while remaining
in a domain of rank two) prior to $t(\varepsilon)$, at $t(\varepsilon)$ it will be in the domain of rank one with probability that tends to one. 
Therefore, as before, for $x \in \overline{D}_6$, $X^{x, \varepsilon}_{t(\varepsilon)}$ tends to 
a uniform distribution on $D_1$, and for $x \in \overline{D}_5$, $X^{x, \varepsilon}_{t(\varepsilon)}$ tends to 
a uniform distribution on $D_2$. However, for $x \in \overline{D}_7 \setminus (\overline{D}_6 \bigcup \overline{D}_5)$, the process will be distributed, in the limit, 
on $D_1 \bigcap D_2$, rather than on $D_3$. To describe
the limiting distribution, we need to consider the process $Y^{x}_t$ corresponding to $k = 7$ and $S = \{3\}$.  Let
$\tau^{x, \varepsilon} = \tau^{x, \varepsilon}(\partial D_6 \bigcup \partial D_5)$ be the first time when $X^{x,\varepsilon}_t$
hits $ \partial D_6 \bigcup \partial D_5$. The stopping time $\tau^x$ is defined similarly, with $Y^{x}_t$ instead of $X^{x,\varepsilon}_t$. If $X^{x,\varepsilon}_t$ reaches 
$\partial D_6$ first, then it will tend to the uniform distribution on $D_1$, and if it reaches $\partial D_5$ first, then it will tend to the uniform
distribution on $D_2$. It will be seen that  $\mathrm{P}(\tau^{x, \varepsilon} \in \partial D_6)$ and $\mathrm{P}(\tau^{x, \varepsilon} \in \partial D_5)$ tend
to the corresponding expressions with $\tau^{x, \varepsilon}$ replaced by $\tau^{x}$. The intuition here is that, even if the process $X^{x,\varepsilon}_t$
reaches $\partial D_3$, there is enough time for it to exit a small neighborhood of $\overline D_3$. Disregarding the time that $X^{x,\varepsilon}_t$ spends in $D_3$,
it is well approximated by $Y^x_t$ until the time it hits $\partial D_6$ or $\partial D_5$. 
Thus, the distribution of $X^{x, \varepsilon}_{t(\varepsilon)}$ tends to 
\[
\mathrm{P}(Y^x_{\tau^{x}} \in \partial D_6) \lambda_1 + \mathrm{P}(Y^x_{\tau^{x}} \in \partial D_5) \lambda_2.
\]

Let us stress that the definition of the stopping time $\tau^x$ here is based on the process $Y^{x}_t$ that is different from the one in (\ref{jijil}).
If $X^{x, \varepsilon}_{t(\varepsilon)}$ starts at $x \in D_3$, the conclusion still holds, with the initial point replaced by an arbitrary point on $\partial D_3$
in order to make sense of $Y^x_t$  and $\tau^x$. (Later, we take the approach where the process $Y^x_t$ is defined on the space where all the points of $\overline{D}_1$ are identified.)

For $x \in D_0 \setminus \overline{D}_7$, the distribution of $X^{x, \varepsilon}_{t(\varepsilon)}$ tends to 
\[
\int_{\partial D_7} ( \mathrm{P}(Y^y_{\tau^{y}} \in \partial D_6) \lambda_1 + \mathrm{P}(Y^y_{\tau^{y}} \in \partial D_5) 
\lambda_2) \pi(x, y) dy,
\]
where $\pi(x, y)$ is still the Poisson kernel for the process in $D_0 \setminus \overline{D}_7$. 

In the case when $\varepsilon^{-2} \ll t(\varepsilon) \ll \varepsilon^{-3}$,  $X^{x, \varepsilon}_{t}$ has enough time to exit a  domain of rank two, but not a small neighborhood of a domain of rank three. Thus, for all $x$,  $X^{x, \varepsilon}_{t(\varepsilon)}$ tends to 
a uniform distribution on $D_1$. 

In the case when $\varepsilon^{-3} \ll t(\varepsilon) \ll \varepsilon^{-4}$, the process has enough time to exit the domain of rank three and will 
visit each of the domains of rank one many times prior to $t(\varepsilon)$. However
$X^{x, \varepsilon}_{t(\varepsilon)}$ still tends to a uniform distribution on $D_1$ as $D_1$ is the `deepest' of all the domains of rank one in the following sense: the time it takes $X^{x, \varepsilon}_t$ to exit a small neighborhood of $D_6 \bigcup D_5 \bigcup D_3$ is much larger for $x \in D_1$ than
for $x \in D_2 \bigcup D_3$. Similarly, $X^{x, \varepsilon}_{t(\varepsilon)}$ still tends to a uniform distribution on $D_1$ for 
 $t(\varepsilon) \gg \varepsilon^{-4}$.

The paper is structured as follows. In Section~\ref{dttpA}, we give a rigorous definition of the process $X^{x,\varepsilon}_t$ with asymmetric behavior on 
the boundaries of the domains $D_k$, $1 \leq k \leq n$. In Section~\ref{dttp}, for each $0 \leq k \leq n$ and appropriate $S$, we define the corresponding process $Y^x_t$ with instantaneous re-distribution and reflection on the boundaries of sub-domains. The main result is formulated in Section~\ref{mremre}. The ingredients necessary
for the proof are developed in Sections~\ref{abstle}-\ref{seconv2}. The proof of the main result is presented in Section~\ref{pmree}. In order to make the exposition more
accessible, we present the proof for the particular example outlined in the Introduction. This example exhibits all the features of the general result, but 
allows us to refer to concrete domains and avoid cumbersome notation.

\section{Processes with asymmetric behavior on the boundaries} \label{dttpA}

We start the discussion with the case of a single domain. Let $D \subset \mathbb{T}^d$ be an open connected domain with infinitely differentiable boundary $\partial D$
and let  $U = \mathbb{T}^d \setminus \overline{D}$. 

The family of processes $X^{x, \varepsilon}_t$, $x \in \mathbb{T}^d$, will be defined in terms of its generator $L^\varepsilon$. Since we expect 
$X^{x, \varepsilon}_t$ to coincide with a Wiener process outside $\partial D$, the generator
coincides with $\frac{1}{2} \Delta$ on a certain class of functions. The domain of the generator, however, should be restricted by certain boundary conditions to account for non-trivial behavior of $X^{x, \varepsilon}_t$ on $\partial D$. We'll use the Hille-Yosida theorem stated here in the form that is convenient for considering closures of linear operators (see \cite{We}).

\begin{theorem}
Let $K$ be a compact space, $C(K)$ be the space of continuous functions on it. The space $C(K)$ is endowed with the supremum norm.  Suppose that a linear operator $A$ on $C(K)$ has the following properties:

(a) The domain $\mathcal{D}(A)$ is dense in $C(K)$;

(b) The constant function  $ \mathbf{1}$ belongs to $\mathcal{D}(A)$ and $A \mathbf{1} = 0$;

(c) The maximum principle: If $S$ is the set of points where a function $f \in \mathcal{D}(A)$ reaches its maximum, then $A f (x) \leq 0$ for at
least one point $x \in S$.

(d) For a dense set $\Psi \subseteq C(K)$, for every $\psi \in \Psi$, and every $\lambda > 0$, there exists a solution $f \in \mathcal{D}(A)$ of the equation $\lambda f - A f = \psi$.

Then the operator $A$ is closeable and its closure $ \overline{A}$ is the infinitesimal generator of a unique semi-group of positivity-preserving operators $T_t$, $t \geq 0$, on $C(K)$ with $T_t \mathbf{1} = \mathbf{1}$, $||T_t|| \leq 1$.
\end{theorem}

The Hille-Yosida theorem will be applied to the space $K = \mathbb{T}^d$. Let us define the linear operator $A^\varepsilon$ in $C(\mathbb{T}^d)$. First, we define its domain. For a function $f \in C(\mathbb{T}^d)$, we denote its restriction to $\overline{U}$ by $f_{\overline{U}}$ and its restriction to $ \overline{D}= 
\mathbb{T}^d \setminus U$ by $f_{\overline{D}}$. For $x \in \partial U$, let $n_U(x)$ be the unit exterior  normal at $x$ (with respect to $U$), and $n_D(x)= -n_U(x)$. 
The domain of $A^\varepsilon$, denoted by $\mathcal{D}(A^\varepsilon)$, consists of all functions $f \in C(\mathbb{T}^d)$ that satisfy the following conditions:

(1) $f_{\overline{U}}$ and $f_{\overline{D}}$ are twice continuously differentiable on $\overline{U}$ and  $ \overline{D}$, respectively. 
%(Both are closed sets with infinitely smooth boundary. The values of the first and second partial derivatives on the two sides of $\partial U$ may be different.);

(2) $\Delta f$ is a continuous function on $\mathbb{T}^d$ (i.e.,  $\Delta f_{\overline{U}}(x) = \Delta f_{\overline{D}}(x)$ for each $ x  \in \partial U$).

(3) $\langle \nabla f_{\overline{U}} (x) , n_U(x) \rangle + \varepsilon^{-1} \langle \nabla f_{\overline{D}} (x) , n_D(x) \rangle = 0,~~x \in \partial D$.
\\
\\
For $f \in \mathcal{D}(A^\varepsilon)$, we define $A^\varepsilon f =  \frac{1}{2} \Delta f$.
\\

Let us check that the conditions of the Hille-Yosida theorem are satisfied.

(a) Consider the set $G$ of functions $g$ that are infinitely differentiable and satisfy 
$\langle \nabla g (x) , n_U(x) \rangle = 0$ for $x \in \partial D$. It is clear that 
$G \subset \mathcal{D}(A)$ and $G$ is dense in $C(\mathbb{T}^d)$.

(b) Clearly $\mathbf{1} \in \mathcal{D}(A)$ and $A^\varepsilon \mathbf{1} = 0$.

(c) If $f$ has a maximum at $x \in U \bigcup D$, it is clear that $ \Delta f (x) \leq 0$. If a maximum is achieved at $x \in \partial D$,
then $\langle \nabla f_{\overline{U}} (x) , n_U(x) \rangle =  \langle \nabla f_{\overline{D}} (x) , n_D(x) \rangle = 0$ (otherwise, one of these two quantities is negative, which
can't happen since there is a maximum at $x$). Therefore, $\nabla f_{\overline{U}} (x) = \nabla f_{\overline{D}} (x) = 0$, which implies that 
$\Delta f_{\overline{U}}(x) = \Delta f_{\overline{D}}(x) \leq 0$.

(d) Let $\Psi$ be the set of infinitely differentiable functions on $\mathbb{T}^d$. It is clear that $\Psi$ is dense in $C(\mathbb{T}^d)$.  The existence of 
a solution $f \in \mathcal{D}(A)$ to the equation 
\[
\lambda f - A f = \psi
\]
can be seen as in  \cite{ET}. The idea of the proof is to consider the functional
\[
F(f) = \int_U (\lambda f^2 - 2 \psi f + \frac{1}{2}  |\nabla f|^2)(x) dx  + \varepsilon^{-1} \int_D (\lambda f^2 - 2 \psi f + \frac{1}{2}  
|\nabla f|^2)(x) dx, ~~~ f \in H^1(\mathbb{T}^d).
\]
It is easy to see that there is a unique $f \in  H^1(\mathbb{T}^d)$ that minimizes $F(f)$. By varying $f$, one then checks that the minimizer 
satisfies the desired differential relation. From standard elliptic theory, it follows that $f$
is sufficiently smooth in $\overline{U}$ and in $\overline{D}$. By varying $f$ in the neighborhood of a boundary point, one then checks that $f$ satisfies the 
required boundary condition. 

Let $ \overline{A^\varepsilon}$ be the closure of $A^\varepsilon$. Let $T^\varepsilon_t$, $t \geq 0$, be the corresponding semi-group on $C(\mathbb{T}^d)$, whose existence is guaranteed by the
Hille-Yosida theorem. By the Riesz-Markov-Kakutani representation theorem, for $x \in \mathbb{T}^d$ there is a measure $P^\varepsilon(t,x,dy)$ on 
$(\mathbb{T}^d, \mathcal{B}(\mathbb{T}^d))$ such that
\[
(T^\varepsilon_t f)(x) = \int_{\mathbb{T}^d} f(y)P^\varepsilon(t,x,dy),~~~f \in C(\mathbb{T}^d).
\]
It is a probability measure since $T^\varepsilon_t \mathbf{1} = \mathbf{1}$. Moreover, it can be easily verified that $P^\varepsilon(t,x,B)$ is a Markov transition function.  Let $X^{x, \varepsilon}_t$, $x \in \mathbb{T}^d$, be the corresponding Markov family. In order to show that a modification with continuous trajectories exists, it is enough to check that $\lim_{t \downarrow 0} P^\varepsilon(t,x,B)/t = 0$ for each closed set $B$ that doesn't contain $x$ (Theorem I.5 of \cite{Mandl}, see also \cite{Dyn}). Let  $f \in \mathcal{D}(A^\varepsilon)$  be a non-negative function that is equal to one on $B$ and whose support doesn't contain $x$. Then
\[
\lim_{t \downarrow 0} \frac{P^\varepsilon(t,x,B)}{t} \leq \lim_{t \downarrow 0} \frac{(T^\varepsilon_t f)(x) - f(x)}{t} = A^\varepsilon f(x) = 0,
\]
as required. Thus $X^{x,\varepsilon}_t$ can be assumed to have continuous trajectories.
\\

Having defined $X^{x,\varepsilon}_t$, let us now discuss some of its basic properties.
Since $ \overline{A^\varepsilon}$ is the infinitesimal generator of the semi-group $T^\varepsilon_t$, we 
have (see Theorem I.1 of \cite{Mandl}), for $f \in \mathcal{D}(A^\varepsilon)$, 
\[
T^\varepsilon_t f - f = \int_0^t T^\varepsilon_s A^\varepsilon f ds,
\]
that is
\[
\mathrm{E} f (X^{x,\varepsilon}_t) - f(x) = \mathrm{E} \int_0^t (A^\varepsilon f) (X^{x,\varepsilon}_s) ds. 
\]
Therefore, since $X^{x,\varepsilon}_t$ is a Markov process with continuous trajectories, for each $x \in \mathbb{T}^d$, the process 
$f (X^{x,\varepsilon}_t) - f(x) - \int_0^t (A^\varepsilon f) (X^{x,\varepsilon}_s) ds$ is a continuous martingale, and, for each stopping time $\tau$ 
with $\mathrm{E} \tau < \infty$, we get
\begin{equation} \label{sttime}
\mathrm{E} f (X^{x,\varepsilon}_\tau) - f(x) = \mathrm{E} \int_0^\tau (A^\varepsilon f) (X^{x,\varepsilon}_s) ds.
\end{equation}
Let $\mu^\varepsilon$ be the measure on $(\mathbb{T}^d, \mathcal{B} (\mathbb{T}^d))$ whose density with respect to the Lebesgue measure $\lambda$ is
\[
p^\varepsilon(x) = \left\{ \begin{array}{ll}
                 1, ~~~~~~ x \in U,\\
                 \varepsilon^{-1}, ~~~ x \in D.
            \end{array}
            \right.
\]
Observe that,  for $f \in  \mathcal{D}(A^\varepsilon)$, 
\[
\int_{\mathbb{T}^d} A^\varepsilon f d \mu^\varepsilon = \frac{1}{2} \int_U  \Delta f_{\overline{U}} d\lambda +  
\frac{\varepsilon^{-1}}{2} \int_D  \Delta f_{\overline{D}} d\lambda = 
\]
\[
\frac{1}{2} \int_{\partial D} (\langle \nabla f_{\overline{U}}  , n_U \rangle + \varepsilon^{-1} \langle \nabla f_{\overline{D}}  , n_D \rangle) d\nu =
 0,
\]
where $\nu$ is the Lebesgue measure on $\partial D$. Since the generator  $ \overline{A^\varepsilon}$  of the process $X^{x,\varepsilon}_t$
is the closure of $A^\varepsilon$, this is enough to conclude (see Theorem 3.37 of \cite{Lig}) that $\mu^\varepsilon$ is invariant for the process, i.e.,
$\mu^\varepsilon(B) = \int_{\mathbb{T}^d} \mathrm{P}(X^{x,\varepsilon}_t \in B) d \mu^\varepsilon(x)$, $B \in \mathcal{B}({\mathbb{T}^d})$. 

 Let us sketch the proof of the fact that the family of processes $X^{x,\varepsilon}_t$, $\varepsilon > 0$, is tight. It is sufficient to check (see \cite{Kt}, Ch. 18) that for each $a, b > 0$ there exists $\delta \in (0,1)$
such that
\begin{equation} \label{ttnns}
\mathrm{P}(\sup_{t \in [0, \delta]} |X^{x,\varepsilon}_t - x| > a) \leq b \delta
\end{equation}
for all $x \in \mathbb{T}^d$ and all $\varepsilon > 0$. 
Let $h(x) = {\rm dist}(x, \partial D)$, and let $r(x) \in \partial D$ be such that ${\rm dist}(x, r(x)) = {\rm dist}(x, \partial D)$. 
The latter function is correctly defined in a small neighborhood of $\partial D$. 
For $y \in \partial D$, 
let $\mathcal{B}^c(y) = \{x: h(x) \leq c, |{\rm dist}(r(x),y)| \leq a/2 \}$. 
Let $\tau(y) = \tau^{c, x, \varepsilon}(y)$ be the first time when the process 
$X^{x,\varepsilon}_t$ starting at $x \in \mathcal{B}^c(y)$ reaches $\partial \mathcal{B}^c(y)$.

Using arguments similar to those employed in the proof of Lemma~\ref{clex}, it is not difficult to show that, for all sufficiently small $c > 0$, 
all $\varepsilon > 0$,
and all $x$ with ${\rm dist}(x, \partial D) \leq c$,
\begin{equation} \label{inrr1}
\mathrm{P}(h(X^{x,\varepsilon}_{\tau(r(x))}) \neq c) \leq b c^4/2.
\end{equation}
Thus, if $c$ is sufficiently small, $X^{x,\varepsilon}_t$ remains,  with probability close to one,  within distance $3a/4$ from the initial point until it reaches a point that is distance $c$ away
from $\partial D$. Since $X^{x,\varepsilon}_t$ coincides with the Brownian motion away from $\partial D$, we also have, for all sufficiently small $c$ and
all $x$ with ${\rm dist}(x, \partial D) \geq c$,
\begin{equation} \label{inrr2}
\mathrm{P}(\sup_{t \in [0,c^4]} |X^{x,\varepsilon}_t - x| > c) \leq b c^4/2.
\end{equation}
Choose $c < a/4$ sufficiently small for (\ref{inrr1}) and (\ref{inrr2}) to hold. 
We obtain (\ref{ttnns}) with $\delta = c^4$ from (\ref{inrr1}) and (\ref{inrr2}) using the strong Markov property of the process.

Let us now generalize the above construction of the process to the case of several (possibly nested) domains inside $\mathbb{T}^d$.
Let $D_1,...,D_n \subset \mathbb{T}^d$ be open connected domains with infinitely differentiable boundaries $\partial D_k$, $k =1,...,n$. The boundaries are assumed to be non-intersecting. Let  $D_0 = \mathbb{T}^d$ .  Define
\[
U_k = D_k \setminus \bigcup_{i: D_i \subset D_k} \overline{D}_i,~~~k =0,...,n.
\]
We assume that there are functions $\varepsilon_k(\varepsilon)$, $k=1,...,n$,  (permeability of $\partial D_k$) taking positive values. 
For a function $f \in C(\mathbb{T}^d)$, we denote its restriction to $\overline{U}_k$ by $f_{\overline{U}_k}$. 
For $x \in \partial U_k$, let $n_{U_k}(x)$ be the unit exterior  normal at $x$ (with respect to $U_k$).

The domain of $A^\varepsilon$, denoted by $\mathcal{D}(A^\varepsilon)$, now consists of all functions $f \in C(\mathbb{T}^d)$ that satisfy the following conditions:

(1) $f_{\overline{U}_k}$ are twice continuously differentiable on $\overline{U}_k$,~~~$k =0,...,n$. 
%(Both are closed sets with infinitely smooth boundary. The values of the first and second partial derivatives on the two sides of $\partial U$ may be different.);

(2) $\Delta f$ is a continuous function on $\mathbb{T}^d$ (i.e.,  $\Delta f_{\overline{U}_k}(x) = \Delta f_{\overline{U}_l}(x)$ for each $ x  \in \partial U_k \bigcap \partial U_l$).

(3) $\langle \nabla f_{\overline{U}_k} (x) , n_{U_k}(x) \rangle + \varepsilon_l^{-1} \langle \nabla f_{\overline{U}_l} (x) , n_{U_l}(x) \rangle = 0,~~x \in \partial U_l \bigcap \partial U_k $, $D_l \subset D_k$.
\\
\\
For $f \in \mathcal{D}(A^\varepsilon)$, we define $A^\varepsilon f =  \frac{1}{2} \Delta f$. As above, it can be checked 
that the conditions of the Hille-Yosida theorem are satisfied, and we can define the process $X^{x,\varepsilon}_t$ with the generator $\overline{A^\varepsilon}$.
The relation (\ref{sttime}) still holds. The invariant measure  $\mu^\varepsilon$ now has the property that its density $p^\varepsilon$ takes a constant
value $p_k^\varepsilon$ on each $U_k$ and $p_k^\varepsilon/p_l^\varepsilon = \varepsilon_{l}$ if $\partial U_l \bigcap \partial U_k \neq \emptyset$ and
 $D_l \subset D_k$. This way, $\mu^\varepsilon$ is defined up to multiplication by a positive constant.  As above, the family 
$X^{x,\varepsilon}_t$, $\varepsilon > 0$, is tight. 
\\

A small modification of the above construction can be used to define
processes with instantaneous reflection at $\partial D_l$ to the interior of $D_l$ (if the process starts outside $\overline{D}_l$, it first reaches $\partial D_l$,
and then continues as a process with reflection to the interior). 
Formally, this corresponds to the situation when $\varepsilon_l = 0$ for some (or all) $l$. 
Such a process, denoted by $Z^x_t$, can be again defined in terms of
its generator $A$: condition (3) is now replaced by 
\\

($3'$) $\langle \nabla f_{\overline{U}_l} (x), n_{U_l}(x) \rangle = 0,~~x \in \partial U_l \bigcap \partial U_k $, $D_l \subset D_k$,
\\
\\ the conditions of the Hille-Yosida theorem are satisfied, 
and the closure of the resulting operator serves as the generator of the process.

\section{Processes with instantaneous re-distribution and reflection on the boundary} \label{dttp}

Let $0 \leq k \leq n$ and $S \subseteq \{1,...,n\}$ be such that $D_l \prec D_k$ for each $l \in S$. 
Let $U = D_k \setminus \bigcup_{l \in S} \overline{D}_l$. We will define the processes $Y^{x}_t$, discussed in the
Introduction, corresponding to the
given values of $k$ and $S$.  If $S$ contains all the indices $l$ such that $D_l \prec D_k$, then  $Y^{x}_t$ will later be identified as the 
limit, as $\varepsilon \downarrow 0$, for the trace of $X^{x,\varepsilon}_t$ (i.e., for the processes obtained from $X^{x,\varepsilon}_t$
by running the clock only when $X^{x,\varepsilon}_t \notin \bigcup_{l \in S} D_l$). Processes with instantaneous re-distribution 
(according to an arbitrary measure) and reflection
were introduced in our earlier work  \cite{FrK1}, \cite{FrK2}. 

 Let $ U = D_k \setminus \bigcup_{l \in S} \overline{D}_l$. Let $U'$ be the metric space obtained from $\overline{U}$ by identifying all points of $\partial D_l$, turning every $\partial D_l$, $l \in S$, into one point $d_l$.  We denote the mapping $\overline{D}_k \rightarrow U'$, where $\overline{D}_l$ gets mapped into $d_l$,
by $h$. 
Clearly, a function $f \in C(U')$ can be viewed as a function on $\overline{U}$ (denoted by $f_{\overline{U}}$) taking constant values on each component $\partial D_l$ of the boundary. For $x \in \partial U$, let $n_U(x)$ be the unit exterior  normal at $x$ (with respect to $U$)

Let  $\nu_l$ be the Lebesgue measure on $\partial D_l$.
The Hille-Yosida theorem will be applied to the space $K = U'$. Let us define the linear operator $L$ in $C(U')$. First we define its domain. It consists of all functions $f \in C(U')$ that satisfy the following conditions:

(1) $f_{\overline{U}}$  is twice continuously differentiable on  $\overline{U}$.

(2) There are constants $g_l$, $l \in S$, such that
\[
\lim_{x \in U, {\rm dist}(x, \partial D_l) \downarrow 0} \Delta f(x) = g_l,~~~l \in S.
\]

(3)  For each $l \in S$,
\begin{equation} \label{intco}
\int_{\partial D_l} \langle \nabla  f_{\overline{U}} , n_U \rangle d \nu_l =0.
\end{equation}

(4) $\langle \nabla  f_{\overline{U}}(x) , n_U(x) \rangle = 0$,~~~$x \in \partial D_k$. 
\\

For $f \in \mathcal{D}(L)$  and $x \in U'$, we define
\[
L f =  \begin{cases} \frac{1}{2} \Delta f(x), & {\rm if }~ x \in U \bigcup \partial D_k, \\ \frac{1}{2} g_l, & {\rm if }~ x = d_l,~~~l \in S. \end{cases}
\]
Let us check that the conditions of the Hille-Yosida theorem are satisfied.

(a) Consider the set $G$ of functions $g$ that are twice continuously  differentiable on $\overline{U}$,
satisfy the relation $\langle \nabla  g_{\overline{U}}(x) , n_U(x) \rangle = 0$ for $x \in \partial D_k$, 
 and have the following property: for each $l \in S$ there is
a set $V_l$ open in $U'$ such that $\partial D_l \subset V_l$ and $g$ is constant on $V_l$. It is clear that $G \subset \mathcal{D}(L)$ and $G$ is dense in $C(U')$.

(b) Clearly $\mathbf{1} \in \mathcal{D}(L)$ and $L \mathbf{1} = 0$.

(c) If $f$ has a maximum at $x \in U \bigcup \partial D_k$, it is clear that $ \Delta f (x) \leq 0$. Now suppose that $f$ has a maximum at $d_l$, $l \in S$.  Note that 
$\langle \nabla f_{\overline{U}}(x) , n_U(x) \rangle$ is identically zero on
$\partial D_l$, since otherwise it would be negative at some points due to (\ref{intco}), which would contradict the fact that $f$ reaches its 
maximum on $\partial D_l$. Then the second derivative of $f_{\overline{U}}$ in the direction of $n$ is non-positive at all points $x \in \partial D_l$. 
Since $f_{\overline{U}}$ is constant on $\partial D_l$, its second derivative in any direction tangential to the boundary is equal to zero. Therefore, 
$\Delta f_{\overline{U}}(x) \leq 0$ for $x \in \partial D_l$, i.e., $L f (d_l) \leq 0$, as required.

(d) Let $\Psi$ be the set of functions $\psi \in C(U')$ that are continuously differentiable on $\overline{U}$. 
It is clear that $\Psi$ is dense in $C(U')$. Let $\widetilde{f} \in C^2 (\overline{U})$ be the solution of the equation
$\lambda \widetilde{f} - \frac{1}{2} \Delta \widetilde{f} = \psi$ in $U$, $\widetilde{f} = 0$ on $\partial D_l$, $l \in S$, 
$\langle \nabla \widetilde{f}(x) , n_U(x) \rangle = 0$, $x \in \partial D_k$. Let $h_l \in C^2(\overline{U})$ be the solution of the equation
\[
\lambda {h}_l(x) - \frac{1}{2} \Delta {h}_l(x) = 0,~~~~x \in U,
\]
\[
h_l(x) = 1,~~x \in \partial D_l;~~~~~h_l(x) = 0,~~x \in \partial U \setminus \partial D_l,
\]
\[
\langle \nabla h_l(x) , n_U(x) \rangle = 0,~~x \in \partial D_k.
\]
Let us look for the solution $f \in \mathcal{D}(L)$ of $\lambda f -  L f = \psi$ in the form $f = \widetilde{f} + \sum_{l \in S}^n c_l h_l$. We get $|S|$ linear equations for $c_l$, $l \in S$. The solution is unique because of the maximum principle. Therefore, the determinant of the system is non-zero, and the solution exists for all the right hand sides.

As before, having verified that the conditions of the  Hille-Yosida theorem are satisfied, we can construct the Markov family $Y^x_t$, $x \in U'$, with continuous trajectories 
whose generator is $ \overline{L}$ (the closure of $L$). 

\section{Formulation of the main result} \label{mremre}
In this section, we will formulate the result on the asymptotic behavior of $X^{x, \varepsilon}_{t(\varepsilon)}$. First, we need to describe the assumptions
on the time scale $t(\varepsilon)$. Suppose that $D_{k_1} \prec D_{k_2} \prec ... \prec D_{k_r}$ and ${\rm rank}(D_{k_r}) = r$. 
In this case, we refer to $C = (D_{k_1},...,D_{k_r})$ as a chain of domains, to $D_{k_1}$ as its first element, and to $D_{k_r}$ as its last element. We'll say that $o_C(\varepsilon) = \varepsilon_{k_1}(\varepsilon) \varepsilon_{k_2}(\varepsilon) ...
\varepsilon_{k_{r}}(\varepsilon)$ is the order of this chain, where we put $\varepsilon_0(\varepsilon) = 1$ for the domain $D_0 = \mathbb{T}^d$ (which is
relevant if $D_0$ is the
last element of the chain). We define the order of a domain $D \in \mathcal{D}$ as
\[
o_{D} (\varepsilon) = \sup_C o_C(\varepsilon),
\]
where the supremum is taken over all chains $C$ whose last element is $D$. Intuitively, $(o_{D}(\varepsilon))^{-1}$ 
is the (order of the) time it takes the process $X^{x, \varepsilon}_{t}$ starting in
$D$ to exit a small neighborhood of this set. (If $D = D_0$, then $(o_C(\varepsilon))^{-1}$ is the supremum over $x$ of the times it takes the process to reach $U_0$).
We make the following assumptions.
\\

{\bf Assumption 1.} For each $k \geq 1$, $\varepsilon_k(\varepsilon) \rightarrow 0$ as $\varepsilon \downarrow 0$. For each pair of chains $C_1$ and $C_2$ with a common last element, either $o_{C_1}(\varepsilon)/o_{C_2}(\varepsilon) \rightarrow 0$ as
$\varepsilon \downarrow 0$ or $o_{C_2}(\varepsilon)/o_{C_1}(\varepsilon) \rightarrow 0$ as
$\varepsilon \downarrow 0$. 

{\bf Assumption 2} $t(\varepsilon) \rightarrow \infty$ as $\varepsilon \downarrow 0$.  For each domain  $D \in \mathcal{D}$, either $o_{D}(\varepsilon)t(\varepsilon) \rightarrow 0$ as
$\varepsilon \downarrow 0$ (in which case $D$ is said to be trapping) or $o_{D}(\varepsilon)t(\varepsilon) \rightarrow \infty$ as
$\varepsilon \downarrow 0$ (in which case $D$ is said to be non-trapping).
%For a domain $D_k$, we define $o_D(\varepsilon) = \max_C o_C(\varepsilon)$, where
%the maximum is taken over 
\\

For $x \in \mathbb{T}^d$, we'll say that $D \in \mathcal{D}$ with ${\rm rank}(D) = r$ is the characteristic domain for $x$ if:

1) $x \in \overline{D}$.

2) Either $D = D_0$ (i.e., $D$ is of maximal rank) or $o_{D}(\varepsilon)t(\varepsilon) \rightarrow 0$ as
$\varepsilon \downarrow 0$.

3) There is no domain with rank lower than $r$ that has properties 1)-2). 
\\

We'll say that $D \in \mathcal{D}$ is a principal
domain if $o_{D}(\varepsilon)/o_{D'}(\varepsilon) \rightarrow 0$ as $\varepsilon \downarrow 0$ whenever  $D', D'' \in \mathcal{D}$ are such that $D \prec D''$ and $D' \prec D''$.
Let $D$ be the characteristic domain for $x \in \mathbb{T}^d$. We'll say that a chain $C = (D_{k_1},..., D_{k_r} = D)$ is admissible if for each $1 \leq i < r$
either $D_{k_i}$ is trapping or it is a principal domain.

Let $C^1,...,C^s$ be the set of all the admissible chains (with the last element $D$ that is the characteristic domain for $x$).
We will see that the limiting distribution
for $X^{x,\varepsilon}_{t(\varepsilon)}$ is a linear combination $c_1\lambda_1 + ... + c_s \lambda_s$ of the uniform distributions $\lambda_1,...,\lambda_s$ concentrated on the first elements of these
chains. The coefficients multiplying the measures $\lambda_i$, $1 \leq i \leq s$, are determined via the following inductive procedure.

Let $\mathcal{D}'$  be the set of all the trapping  domains such that $D' \prec D$ for $D' \in \mathcal{D}'$. 
Similarly, let $\mathcal{D}''$  be the set of all the non-trapping  domains such that $D'' \prec D$ for $D'' \in \mathcal{D}''$. 
If $\mathcal{D}'$ is empty, then there is only one admissible chain, and the limiting distribution is concentrated on the
first element of this chain.

Next, we describe the coefficients $c_i$ under the assumption that $\mathcal{D}'$ is non-empty, however, each $D' \in \mathcal{D}'$ does
not contain trapping sub-domains. 
 Let $Y^{h(x)}_t$ be the process in the space $U'$  corresponding to domain $D$ and the collection of subdomains 
$\mathcal{D}''$ (as in Section~\ref{dttp}).   Let $\tau$ be the first time this process reaches $\bigcup_{D' \in \mathcal{D}'} \partial D'$. 
Let $\pi$ be the measure induced by $Y^{h(x)}_\tau$ on $\bigcup_{D' \in \mathcal{D}'} \partial D'$. The number of admissible chains is equal to 
the number of elements in $\mathcal{D}'$ (the $i$-th chain has some $D'_i \in \mathcal{D}'$ as its next-to-last element). We claim
that $c_i = \pi(\partial D'_i)$. 

Finally, assume that we know how to determine $c_i = c_i(x)$ for each $x \in \overline{D'}$, where $D' \in \mathcal{D}'$, in the case when $\mathcal{D}'$ is non-empty.
Define the process $Y^{h(x)}_t$, the stopping time $\tau$ and the measure $\pi$ as above.
We claim that
\[
c_i(x) = \int_{\bigcup_{D' \in \mathcal{D}'} \partial D'} c_i(y) d\pi(y),
\]
where the integrand in the right hand side is defined by our inductive assumption. 

\begin{theorem} \label{mteor} Suppose that Assumptions 1 and 2 are satisfied. 
Let $C^1,...,C^s$ be the set of all the admissible chains (with the last element $D$ that is the characteristic domain for $x$). The limiting distribution
for $X^{x,\varepsilon}_{t(\varepsilon)}$ is a linear combination $c_1\lambda_1 + ... + c_s \lambda_s$ of the uniform distributions $\lambda_1,...,\lambda_s$ concentrated on the first elements of these
chains. The coefficients multiplying the measures $\lambda_i$, $1 \leq i \leq s$, are determined via the inductive procedure described
above.
\end{theorem} 

\section{A lemma on the convergence of processes} \label{abstle}

In this section, we prove a lemma that will be useful for establishing the convergence of the trace of $X^{x,\varepsilon}_t$
to the process  $Y^{x}_t$ defined in Section~\ref{dttp}. To simplify the discussion, let us assume that $k = 0$, i.e., $D_k = D_0 = \mathbb{T}^d$, and 
$S$ contains all the indices $l$ such that $D_l \prec D_0$. Recall that 
$U = D_0 \setminus \bigcup_{l \in S} \overline{D}_l$ and $h: \mathbb{T}^d \rightarrow U'$ is the mapping defined by $h(x) = x$ for $x \in U$ and $h(x) = d_l$ for $x \in \overline{D}_l$.  
%\br Maybe distinguish between $U_0 = U$ and $U_k$ at some point. \be
For each $t \geq 0$, define the stopping time
\[
s(t) = \inf(s: \lambda(u: u \leq s, X^{x, \varepsilon}_u \in \overline{U}) = t)),
\]
where $\lambda$ is the Lebesgue measure on the real line, and let
\[
Y^{x, \varepsilon}_t = X^{x, \varepsilon}_{s(t)}.
\]
Thus $Y^{x, \varepsilon}_t $ is a left-continuous process with values in $\overline{U}$, which also can be viewed as a  continuous $U'$-valued process. It can be obtained from $X^{x, \varepsilon}_t $ by running the clock only when $X^{x, \varepsilon}_t $ is in $\overline{U}$.

Note that while convergence of $Y^{x, \varepsilon}_t$ to Markov processes on $U'$ as $\varepsilon \downarrow 0$ will be established, the processes $Y^{x, \varepsilon}_t$
need not be Markov for fixed $\varepsilon > 0$. The main point of the next lemma is that,  in order to demonstrate the convergence of 
$Y^{x, \varepsilon}_t$ to a limiting process, it is sufficient to check that for small $\varepsilon$ the processes
nearly satisfy the relation (\ref{mprobaa}), which is similar to the martingale problem but with the ordinary expectation rather than the conditional 
expectation. A similar lemma (in the situation that did not involve the time change, however) was used in \cite{FW}, Ch 8.

\begin{lemma} \label{flei}  Let $Y^y_t$, $y \in U'$, be a Markov family on $U'$ with continuous trajectories whose semigroup 
$T_t$, $t \geq 0$,  preserves the space $C(U')$. Let $L: D(L) \rightarrow C(U')$ denote the infinitesimal generator of this family,  where $D(L)$ is the domain of the generator.
Let $\Psi$ be a dense linear subspace  of $C(U')$ and $D$ be a linear
subspace  of $D(L)$, and
suppose that $\Psi$ and $D$ have the following properties:

(1)  There is   $\lambda > 0$ such that  for  each $f \in \Psi$ the equation $\lambda F -  L F = f$ has a solution~$F \in D$. 

(2) For each $t > 0$, each $f \in D$, 
  \begin{equation} \label{mprobaa}
	\lim_{\varepsilon \downarrow 0} \mathrm{E} (f( Y^{x, \varepsilon}_t ) - f (  Y^{x, \varepsilon}_0) - \int_0^t L f (  Y^{x, \varepsilon}_u) du) =  0,
	\end{equation}
uniformly in $x \in \mathbb{T}^d$. 

%Suppose that the family of measures on $C([0, \infty), U')$ induced by the processes $Y^{x, \varepsilon}_t$, $\varepsilon > 0$, is tight for each $x \in \mathbb{T}^d$. 
 
Then, for each $x \in \mathbb{T}^d$, the measures induced by the processes $Y^{x, \varepsilon}_t$ converge weakly, as $\varepsilon \downarrow 0$, to
the measure induced by the process $Y^{h(x)}_t$. 
\end{lemma}
\proof Fix $x \in \mathbb{T}^d$. Observe that the family of measures on $C([0, \infty), U')$ induced by the processes $Y^{x, \varepsilon}_t$, $\varepsilon > 0$, is tight since the processes coincide with a Brownian motion on $U$ (and $U' \setminus U$ consists of a finite set of points). Therefore,
we can find a process $\bar{Y}^x_t$ with continuous trajectories and a sequence
$\varepsilon_n \downarrow 0$ such that  $Y^{x, \varepsilon_n}_t$ converge to $\bar{Y}^x_t$ in distribution as $n \rightarrow \infty$. The desired result will immediately follow if we demonstrate that the distribution of $\bar{Y}^x_t$ 
coincides with the distribution of $Y^{h(x)}_t$ (and thus does not depend
on the choice of the sequence $\varepsilon_n$). We will show that $\bar{Y}^x_t$ is a solution of the martingale problem for $(L|_D, h(x))$, i.e.,  for each $t_2 > t_1 \geq 0$ and $f \in D$,
  \begin{equation} \label{rmpro}
\mathrm{E} (f( \bar{Y}^{x}_{t_2} ) - f (  \bar{Y}^{x}_{t_1}) - \int_{t_1}^{t_2} L f (  \bar{Y}^{x}_u) du| {\mathcal{F}}_{t_1}^{\bar{Y}^x}) =  0,~~~~~\bar{Y}^x_0 = h(x). 
	\end{equation}
First, however, let us discuss the uniqueness for solutions of the martingale problem. 	
We claim that:

(a) $D$ is dense in $C(U')$.

(b) ${\rm Range}(\lambda - L|_D)$ is dense in $C(U')$.

(c) For each pair of measures $\mu_1$, $\mu_2$ on $U'$, the equality $\int_{U'} f d\mu_1 = \int_{U'} f d \mu_2$ for all $f \in C(U')$ implies that
$\mu_1 = \mu_2$.

 To demonstrate (a), 
take an arbitrary $\delta > 0$ and $F_0 \in D(L)$. Let $g_0 = \lambda F_0 - L F_0$, and take $g' \in \Psi$ such that $\|g' - g_0\| \leq \lambda \delta$. Let $F' \in D$
be such that $\lambda F' - L F' = g'$. Then, since $L$ is the generator of a strongly continuous semigroup on $C(U')$, from the Hille-Yosida theorem  
it follows that $\|F' - F_0\| \leq \|g' - g_0\|/\lambda \leq \delta$. This implies (a) since  $D(L)$ is dense in $C(U')$. Note that (b) follows from the existence of a solution $F \in D$ to $\lambda F -  L F = f \in \Psi$ and the density of $\Psi$, while (c) is obvious.  The validity of (a)-(c) is enough to conclude that the distribution on $C([0, \infty), U')$ of a process with continuous paths satisfying (\ref{rmpro}) is uniquely determined (Theorem 4.1, Chapter 4 in \cite{EK86}).

Note that (\ref{rmpro}) is satisfied if $\bar{Y}^x_t$ is replaced by $Y^{h(x)}_t$ since $D \subseteq D(L)$ and $L$ the the generator of the family $Y^y_t$, $y \in U'$. 
Therefore, $\bar{Y}^x_t$ and $Y^{h(x)}_t$ have the same distribution if (\ref{rmpro}) holds. It remains to prove (\ref{rmpro}).

Note that $\bar{Y}^x_t$ is a solution of the martingale problem for $(L|_D, h(x))$ if and only if 
\[
\mathrm{E} \left( (\prod_{i=1}^k  g_i(\bar{Y}^x_{u_i})) (f( \bar{Y}^{x}_{t_2} ) - f (  \bar{Y}^{x}_{t_1}) - \int_{t_1}^{t_2} L f (  \bar{Y}^{x}_u) du) \right) =  0,~~~~~\bar{Y}^x_0 = h(x),
\]
whenever $f \in D$, $0 \leq u_1 < ... < u_k \leq t_1$, and $g_1,...,g_k \in C(U')$. Since $Y^{x, \varepsilon_n}_t = h(X^{x, \varepsilon_n}_{s(t)})$ converge to $\bar{Y}^x_t$ in distribution, we have
\[
\mathrm{E} \left( (\prod_{i=1}^k  g_i(\bar{Y}^x_{u_i})) (f( \bar{Y}^{x}_{t_2} ) - f (  \bar{Y}^{x}_{t_1}) - \int_{t_1}^{t_2} L f (  \bar{Y}^{x}_u) du) \right) =
\]
\[
\lim_{n \rightarrow \infty} \mathrm{E} \left( (\prod_{i=1}^k  g_i(h(X^{x,\varepsilon_n}_{s(u_i)}))) (f( h(X^{x,\varepsilon_n}_{s(t_2)}) ) - 
f (  h(X^{x,\varepsilon_n}_{s(t_1)})) - \int_{t_1}^{t_2} L f (  h(X^{x,\varepsilon_n}_{s(u)})) du) \right) = 
\]
\[
\lim_{n \rightarrow \infty} \mathrm{E} \left( (\prod_{i=1}^k  g_i(h(X^{x,\varepsilon_n}_{s(u_i)}))) \mathrm{E}(f( h(X^{x,\varepsilon_n}_{s(t_2)}) ) - f (  h(X^{x,\varepsilon_n}_{s(t_1)})) - \int_{t_1}^{t_2} L f (  h(X^{x,\varepsilon_n}_{s(u)})) du| \mathcal{F}_{s(t_1)}^{X^{x,\varepsilon_n}}) \right).
\]
By the strong Markov property of the family $X^{x, \varepsilon_n}_t$, 
\[
\mathrm{E}(f( h(X^{x,\varepsilon_n}_{s(t_2)}) ) - f (  h(X^{x,\varepsilon_n}_{s(t_1)})) - \int_{t_1}^{t_2} L f (  h(X^{x,\varepsilon_n}_{s(u)})) du| 
\mathcal{F}_{s(t_1)}^{X^{x,\varepsilon_n}}) = 
\]
\[
\mathrm{E}(f( h(X^{x',\varepsilon_n}_{s(t_2-t_1)}) ) - f (  h(X^{x',\varepsilon_n}_{0})) - \int_{0}^{t_2-t_1} L f (  h(X^{x',\varepsilon_n}_{s(u)})) du) 
|_{x' = X^{x,\varepsilon_n}_{s(t_1)}}~, 
\]
which tends to zero in distribution, as follows from (\ref{mprobaa}).
Therefore, using the boundedness of $f$, $Lf$, and $g_1,...,g_k$, we conclude that 
\[
\mathrm{E} \left( (\prod_{i=1}^k  g_i(\bar{Y}^x_{u_i})) (f( \bar{Y}^{x}_{t_2} ) - f (  \bar{Y}^{x}_{t_1}) - \int_{t_1}^{t_2} L f (  \bar{Y}^{x}_u) du) \right) = 0.
\]
Finally, $\bar{Y}^x_0 = h(x)$ since $Y^{x, \varepsilon_n}_0 = h(X^{x,\varepsilon_n}_0) = h(x)$ for all $n$. \qed
\\

In order to deal with the case when $k \neq 0$, i.e., $D_k \neq \mathbb{T}^d$, we need to understand the behavior of the process $X^{x,\varepsilon}_t$
near the boundary of $\partial D_k$. The following lemma will be useful in  proving the convergence of  $X^{x,\varepsilon}_t$ to the reflected Brownian motion
in the case when $D_k$ does not contain sub-domains (is a domain of rank one). 
This lemma is similar to Lemma~\ref{flei}, but now there is no time change. 
\begin{lemma}  \label{ntch} Let $Z^x_t$, $x \in \mathbb{T}^d$, be a Markov family on $\mathbb{T}^d$ with continuous trajectories whose semigroup 
$T_t$, $t \geq 0$,  preserves the space $C(\mathbb{T}^d)$. Let $L: D(L) \rightarrow C(\mathbb{T}^d)$ denote the infinitesimal generator of this family,  where $D(L)$ is the domain of the generator.
Let $\Psi$ be a dense linear subspace  of $C(\mathbb{T}^d)$ and $D$ be a linear
subspace  of $D(L)$, and
suppose that $\Psi$ and $D$ have the following properties:

(1)  There is   $\lambda > 0$ such that  for  each $f \in \Psi$ the equation $\lambda F -  LF = f$ has a solution~$F \in D$. 

(2) For each $t > 0$, each $f \in D$, 
  \[
	\lim_{\varepsilon \downarrow 0} \mathrm{E} (f( X^{x, \varepsilon}_t ) - f (  X^{x, \varepsilon}_0) - \int_0^t L f (  X^{x, \varepsilon}_u) du) =  0,
	\]
uniformly in $x \in \mathbb{T}^d$. 

%Suppose that the family of measures on $C([0, \infty), U')$ induced by the processes $Y^{x, \varepsilon}_t$, $\varepsilon > 0$, is tight for each $x \in \mathbb{T}^d$. 
 
Then, for each $x \in \mathbb{T}^d$, the measures induced by the processes $X^{x, \varepsilon}_t$ converge weakly, as $\varepsilon \downarrow 0$, to
the measure induced by the process $Z^x_t$. 
\end{lemma}

This lemma is proved in the same way as Lemma~\ref{flei}. The only difference is that, while there it was obvious that the family of measures on $C([0, \infty), U')$ 
induced by the processes $Y^{x, \varepsilon}_t$, $\varepsilon > 0$, was tight, to claim tightness for the family of measures on $C([0, \infty), \mathbb{T}^d)$
induced by the processes  $X^{x, \varepsilon}_t$, $\varepsilon > 0$, we now need to refer to Section~\ref{dttpA}.

\section{Behavior of the process near the boundary of a trapping domain} \label{bounarybeh}
Consider the case of a single trapping domain $D_1 = D \subset D_0 = \mathbb{T}^d$. Assume that $\varepsilon_1 (\varepsilon)  =\varepsilon$. 
%In this section, we prove several lemmas concerning the behavior of $X^{x, \varepsilon}_t$ near $\partial D$. 
Let $S_r = \{x \in \overline{U}: {\rm dist}(x, \partial D) = r\}$ for $r \geq 0$,  $S_r = \{x \in D: {\rm dist}(x, \partial D) = -r\}$ for $r < 0$. These are smooth surfaces if $r$ is sufficiently small.  For $r > 0$, let
\[
\Gamma_r = \{x \in \mathbb{T}^d: {\rm dist}(x, \partial D) \leq r \},~~\Gamma^+_r = \{x \in \overline{U}: {\rm dist}(x, \partial D) \leq r \}.
\]
% Thus \br $  \tau^{x, \varepsilon} = \tau^{x, \varepsilon}(\partial D)$. We also define $\sigma^{x, \varepsilon} = \tau^{x, %\varepsilon}(S_{\sqrt{\varepsilon}})$.
%\be 
For $B \subset \mathbb{T}^d$, let $  \tau^{x, \varepsilon}(B) =  \inf \{t\geq 0: X^{x,\varepsilon}_t \in B \}$.

We will see that if the process starts at $x \in \partial D$, then, with probability close to one, 
it exits $\Gamma_{{\varepsilon}^\alpha}$ in a location that is close to $x$. Let $h(x) = {\rm dist}(x, \partial D)$, and let $r(x) \in \partial D$ be such that ${\rm dist}(x, r(x)) = {\rm dist}(x, \partial D)$. 
The latter function is correctly defined in a small neighborhood of $\partial D$.
Let $0 < \beta < \alpha < 1 $ and 
 $\mathcal{G}_\varepsilon(x) = \{y: h(y) \leq \varepsilon^\alpha, |r(y) - x| \leq \varepsilon^\beta \}$. 
\begin{lemma} \label{clex}
For each $0 < \beta < \alpha < 1 $,
\[
\mathrm{P} (X^{x, \varepsilon}_{\tau^{x, \varepsilon}(\partial \mathcal{G}_\varepsilon(x))}  \in S_{-{\varepsilon^\alpha}} \bigcup S_{{\varepsilon^\alpha}}) \rightarrow 1
\]
as $\varepsilon \downarrow 0$ uniformly in $x \in \partial D$. 
\end{lemma}
\proof For $y \in \mathcal{G}_\varepsilon(x)$, define $f(y) =  \varepsilon^{-2\beta} ({\rm dist}(r(y), x))^2 - c  h^2(y)) $. Here the constant $c$ is chosen so
large that $\Delta f(y) \leq 0$ for $y \in \mathcal{G}_\varepsilon(x)$.  We extend $f$ to $\mathbb{T}^d$ so that $f \in \mathcal{D}(A^\varepsilon)$ and apply
 (\ref{sttime}) with $\tau = \tau^{x, \varepsilon}(\partial \mathcal{G}_\varepsilon(x) )$. Thus
\[
\mathrm{E} f ( X^{x, \varepsilon}_{\tau^{x, \varepsilon}(\partial \mathcal{G}_\varepsilon(x))}) \leq 0. 
\]
Observe that $f$ is bounded from below on $\partial \mathcal{G}_\varepsilon(x) $ by $- c \varepsilon^{2(\alpha - \beta)}$. Therefore,
\[
\mathrm{E} \max(0, f ( X^{x, \varepsilon}_{\tau^{x, \varepsilon}(\partial \mathcal{G}_\varepsilon(x))}) ) \leq  c \varepsilon^{2(\alpha - \beta)}.
\]
Since $f(y) \geq 1/2$ on $\partial \mathcal{G}_\varepsilon(x) \setminus (S_{-{\varepsilon^\alpha}} \bigcup S_{{\varepsilon^\alpha}})$ for all sufficiently small $\varepsilon$, we conclude that
\[
\mathrm{P} (X^{x, \varepsilon}_{\tau^{x, \varepsilon}(\partial \mathcal{G}_\varepsilon(x))}  \notin S_{-{\varepsilon^\alpha}} \bigcup S_{{\varepsilon^\alpha}})
\leq 2 c \varepsilon^{2(\alpha - \beta)} \rightarrow 0,
\]
which gives the desired result.
\qed
\\

 The next lemma provides an estimate on the time it takes the process $X^{x, \varepsilon}_t$ starting at $x \in \partial D$ to exit $\Gamma_{{\varepsilon^\alpha}}$.
\begin{lemma} \label{timeexit}  For each $\alpha \in (0,1)$, there is a constant $c > 0$ such that
\[
\sup_{x \in \partial D} \mathrm{E} \tau^{x, \varepsilon}(S_{-{\varepsilon^\alpha}} \bigcup S_{{\varepsilon^\alpha}}) \leq c \varepsilon^{2\alpha}
\]
for all sufficiently small $\varepsilon$. 
\end{lemma}
\proof  Recall the definition of the operator $A^\varepsilon$  from Section~\ref{dttpA}. Since $\partial D$ is smooth, for sufficiently small $r > 0$, there
exists $f \in \mathcal{D}(A^\varepsilon)$ satisfying $f(x) = ({\rm dist}(x, \partial D))^2$ when ${\rm dist}(x, \partial D) \leq r$. The lemma immediately follows
from (\ref{sttime}) with $\tau = \tau^{x, \varepsilon}(S_{-{\varepsilon^\alpha}} \bigcup S_{{\varepsilon^\alpha}})$. 
\qed
\\

We can control the probability with which the process exits $\Gamma_{{\varepsilon}^\alpha}$ through $S_{{\varepsilon}^\alpha}$. 
\begin{lemma} \label{dirofexit}
 For each $\alpha \in (0,1)$,
\[
\lim_{\varepsilon \downarrow 0} \varepsilon^{-1} \mathrm{P}  
(X^{x,\varepsilon}_{\tau^{x, \varepsilon}(S_{-{\varepsilon}^\alpha} \bigcup S_{{\varepsilon^\alpha}})} \in S_{{\varepsilon^\alpha}}) = 1
\]
uniformly in $x \in \partial D$.
\end{lemma}
\proof  For sufficiently small $r > 0$, define the function $f$ in the $\Gamma_r$:  $f(x) = 0$ for $x \in \partial D$;
$f(x) = \varepsilon ({\rm dist}(x, \partial D) + g(x) ({\rm dist}(x, \partial D))^2)$ for $x \in D \bigcap \Gamma_r$; 
$f(x) = - {\rm dist}(x, \partial D) + g(x) ({\rm dist}(x, \partial D))^2$ for $x \in U \bigcap \Gamma_r$. 
 The function $g \in C^2(\Gamma_r)$  can be chosen in such a way that $\lim_{y \rightarrow x, y \in D} \Delta f (y) = 
\lim_{y \rightarrow x, y \in U} \Delta f (y) = 0$
for each $x \in \partial D$ ($g$ can be first defined on $\partial D$ and assumed to be constant on each segment perpendicular to $\partial D$).  We continue $f$ outside $\Gamma_r$ so that $f \in \mathcal{D}(A^\varepsilon)$.
Applying  (\ref{sttime}) with $\tau = \tau^{x, \varepsilon}(S_{-{\varepsilon}^\alpha} \bigcup S_{{\varepsilon}^\alpha})$, we obtain
\[
\mathrm{E} f (X^{x,\varepsilon}_\tau) = \frac{1}{2} \mathrm{E} \int_0^\tau   \Delta f (X^{x,\varepsilon}_s) ds.
\]
Therefore, using Lemma~\ref{timeexit} to estimate the integral in the right hand side, we obtain
\[
\mathrm{P}  (X^{x,\varepsilon}_{\tau} \in S_{{\varepsilon}^\alpha}) (1 + o(1)) - 
\varepsilon \mathrm{P}  (X^{x,\varepsilon}_{\tau} \in S_{-{\varepsilon}^\alpha}) (1 + o(1))  = O(\varepsilon^{2\alpha}).
\]
This shows that $\mathrm{P}  (X^{x,\varepsilon}_{\tau} \in S_{{\varepsilon}^\alpha}) \rightarrow 0$ and, therefore, $ \mathrm{P}  (X^{x,\varepsilon}_{\tau} \in S_{-{\varepsilon}^\alpha}) \rightarrow 1$. The same formula now immediately implies the statement of the lemma under the additional condition that $\alpha > 1/2$.
For $\alpha \in (0, 1/2]$, we can use the validity of the lemma for $\alpha' = 3/4$, and the strong Markov property of the process. (In order to reach 
$S_{-{\varepsilon}^\alpha} \bigcup S_{{\varepsilon^\alpha}}$, the process must first reach $S_{-{\varepsilon}^{\alpha'}} \bigcup S_{{\varepsilon^{\alpha'}}}$, while upon reaching
the latter, it either returns to $\partial D$ or proceeds to $S_{-{\varepsilon}^\alpha} \bigcup S_{{\varepsilon^\alpha}}$. The probability
of the latter event, given a starting point in $S_{-{\varepsilon}^{\alpha'}} \bigcup S_{{\varepsilon^{\alpha'}}}$, is  
asymptotically equivalent to $\varepsilon^{\alpha - \alpha'}$, uniformly in the starting point, since the process coincides with the Brownian motion 
outside $\partial D$.) \qed
\\
% We can also estimate the time the process $X^{x, \varepsilon}_t$ spends in  $\Gamma_{{\varepsilon^\alpha}}^+$ prior to exiting $\Gamma_{{\varepsilon^\alpha}}$.

Next, we estimate the time spent by the process in $\Gamma^+_{\varepsilon^\alpha}$ prior to reaching $S_{\varepsilon^\alpha}$.
\begin{lemma} \label{stwo}
For each $\alpha \in (0,1)$, there is $c >0$ such that
\begin{equation} \label{eqmm}
 \mathrm{E} \lambda(t: X^{x, \varepsilon}_t \in \Gamma^+_{\varepsilon^\alpha}, ~ 0 \leq t \leq \tau^{x, \varepsilon}(S_{\varepsilon^\alpha}))  
\leq c\varepsilon^{2 \alpha}.
\end{equation}
for all $x \in \overline{D}$, where $\lambda$ is the Lebesgue measure on the real line. 
\end{lemma}
\proof  Consider a function $f$ that satisfies: $f(x) = ({\rm dist}(x, \partial D))^2$ when ${\rm dist}(x, \partial D) \leq r$, $x \in U$; $f \in C^2(U)$;
$f (x) = 0$, $x \in \overline{D}$. This function does not belong to $\mathcal{D}(A^\varepsilon) $  since $\Delta f$ is not continuous. However, there exist
functions $f_n \in \mathcal{D}(A^\varepsilon)$ such that $\Delta f_n$ are uniformly bounded; $\Delta f_n (x) \rightarrow \Delta f(x)$ for all $x \notin \partial D$;
$f_n(x) \rightarrow f(x)$ for all $x$. Therefore, since (\ref{sttime}), with $\tau = \tau^{x, \varepsilon}(S_{-{\varepsilon^\alpha}} \bigcup S_{{\varepsilon^\alpha}})$,
is applicable to $f_n$, it is also applicable to $f$. Therefore, by Lemma~\ref{dirofexit}, for each $\alpha \in (0,1)$, there is a constant $c > 0$ such that
\[
\sup_{x \in \partial D}   \mathrm{E} \lambda(t: X^{x, \varepsilon}_t \in \Gamma^+_{{\varepsilon}^\alpha}, ~ 0 \leq t \leq \tau^{x, \varepsilon}(S_{-{\varepsilon^\alpha}} \bigcup S_{{\varepsilon^\alpha}})) \leq c \varepsilon^{2\alpha+1},
\]
where $\lambda$ is the Lebesgue measure on the real line. 

Let us now return to the proof of (\ref{eqmm}). Without loss of generality, we may assume that $x \in \partial D$.  
Let $\sigma^{x, \varepsilon}_0 = 0$,  $\tau^{x, \varepsilon}_n = \inf(t \geq \sigma^{x, \varepsilon}_{n-1}: X^{x, \varepsilon}_t 
\in S_{-{\varepsilon}^{\alpha}} \bigcup S_{{\varepsilon^{\alpha}}})$, $n \geq 1$,  while
$\sigma^{x, \varepsilon}_n = \inf(t \geq \tau^{x, \varepsilon}_n: X^{x, \varepsilon}_t \in \partial D)$, 
$n \geq 1$. Let $N^{x,\varepsilon}  = \min(n: X^{x, \varepsilon}_{\tau_n} \in S_{{\varepsilon}^\alpha})$. Thus $N^{x,\varepsilon}$ is the number of excursions prior to
reaching $ S_{{\varepsilon}^\alpha}$. By Lemma~\ref{dirofexit},
\[
 \lim_{\varepsilon \downarrow 0} ( \varepsilon \mathrm{E} N^{x,\varepsilon})  = 1.
\]
Combining this with the above bound on the expected contribution from one excursion, 
we obtain the desired result.
\qed
\\

\begin{lemma} \label{locexit} For each $f \in C^2(\overline{D}) \bigcap C(\mathbb{T}^d)$, $\alpha \in (0,1)$,
\[
\sup_{x \in \partial D} | \frac{\mathrm{E}f(X^{x, \varepsilon}_{\tau^{x, \varepsilon}(S_{-{\varepsilon}^\alpha} \bigcup S_{{\varepsilon}^\alpha})}) - 
f(x)}{\varepsilon^\alpha} - 
 \langle \nabla f_{\overline{D}} (x) , n_D(x) \rangle|  \rightarrow 0,~~ as~~\varepsilon \downarrow 0.
\]
\end{lemma}
\proof 
%Let $x = x_0 \in \partial D$ be the initial point for the process  $X^{x, \varepsilon}_t$. We'll treat $x_0$ as fixed, but all the estimates and limits below 
%are easily seen to be uniform in $x_0$. 
Let $h(x) = {\rm dist}(x, \partial D)$, and let $r(x) \in \partial D$ be such that ${\rm dist}(x, r(x)) = {\rm dist}(x, \partial D)$. 
The latter function is correctly defined in a small neighborhood of $\partial D$.
%Fix $\beta \in (\frac{1}{2}, \alpha)$. 
% Let $\mathcal{G}_\varepsilon = \{x: h(x) \leq \varepsilon^\alpha, |r(x) - x_0| \leq \varepsilon^\beta \}$. 
Let us define $\tilde{f}: \Gamma_\varepsilon \rightarrow \mathbb{R}$ by putting $\tilde{f}(x) = f(r(x))$. We extend $\tilde{f}$ to $\mathbb{T}^d$ as
a function from  $\mathcal{D}(A^\varepsilon)$.  Applying  (\ref{sttime}) 
with $\tau = \tau^{x, \varepsilon}(S_{-{\varepsilon}^\alpha} \bigcup S_{{\varepsilon}^\alpha})$ to $\tilde{f}$, we obtain
\[
\mathrm{E} \tilde{f} ((X^{x,\varepsilon}_\tau) - \tilde{f}(x)) = \frac{1}{2} \mathrm{E} \int_0^\tau    \Delta  \tilde{f} (X^{x,\varepsilon}_s) ds.
\]
By Lemma~\ref{timeexit}, the absolute value of the right hand side is bounded from above by $c \varepsilon^{2\alpha}$. 
Therefore,
\[
\frac{\mathrm{E}f(X^{x, \varepsilon}_{\tau}) - 
f(x)}{\varepsilon^\alpha} - 
 \langle \nabla f_{\overline{D}} (x) , n_D(x) \rangle  =
\]
\[
 \mathrm{E}\frac{f(X^{x, \varepsilon}_{\tau}) - 
\tilde{f}(X^{x, \varepsilon}_{\tau})}{\varepsilon^\alpha} - 
 \langle \nabla f_{\overline{D}} (x) , n_D(x) \rangle + O(\varepsilon^\alpha) =
\]
\[
\mathrm{E}\frac{f(X^{x, \varepsilon}_{\tau}) - 
f(r(X^{x, \varepsilon}_{\tau}))}{\varepsilon^\alpha} - 
 \langle \nabla f_{\overline{D}} (x) , n_D(x) \rangle + O(\varepsilon^\alpha).
\]
The right hand side tends to zero uniformly in $x \in \partial D$, as follows from Lemma~\ref{dirofexit} and the fact that $X^{x, \varepsilon}_{\tau} \rightarrow x$ 
in probability (by Lemma~\ref{clex}). 
%
%Without loss of generality, we can assume that $f(x_0) = 0$. Let $\tilde{f}$ be the linear function whose derivatives in the directions tangential to $\partial D$ at %$x_0$ coinside with
%those of $f$, and the normal derivative is zero. We assume that $\tilde{f}(x_0) = f(x_0) = 0$. Let $\tilde{f}_1$ be the restriction
%of $\tilde{f}$ to $\overline{D}\bigcap \mathcal{G}_\varepsilon$. We would like to us extend  $\tilde{f}_1$ to a function that belongs to $\mathcal{D}(A^\varepsilon)$. %First, we
%define $\tilde{f}_2$ on $\overline{U}\bigcap \mathcal{G}_\varepsilon$ via
%\[
%\tilde{f}_2 (x) = \tilde{f}(x) + \varepsilon^{-1}(\langle \nabla \tilde{f} (r(x)) , n_D(r(x)) \rangle h(x) + g(r(x))h^2(x)), 
%\]
%\br where a smooth function $g$ is chosen in such a way that $\Delta \tilde{f}_2 (x) = 0$ for $x \in \partial D \bigcap \mathcal{G}_\varepsilon$. \be
%Define $\hat{f}$ to be equal to $\tilde{f}_1$ in $\overline{D}\bigcap \mathcal{G}_\varepsilon$, $\tilde{f}_2$ in $\overline{U}\bigcap \mathcal{G}_\varepsilon$,
%and extend it as an arbitrary element of $\mathcal{D}(A^\varepsilon)$ outside $\mathcal{G}_\varepsilon$. Applying  (\ref{sttime}) 
%with $\tau = \tau^{x, \varepsilon}(\mathcal{G}_\varepsilon)$ to $\hat{f}$, we obtain
%\begin{equation} \label{innn4}
%\mathrm{E} \hat{f} (X^{x_0,\varepsilon}_\tau) = \mathrm{E} \int_0^\tau (A^\varepsilon \hat{f}) (X^{x_0,\varepsilon}_s) ds.
%\end{equation}
%
%%
\qed
\\

\section{Exit from a neighborhood of a trapping domain} \label{exittd}
Consider the case of a single trapping domain $D$. First, we estimate the time it takes the process starting at
$x \in \overline{D}$ to exit a small neighborhood of $\overline{D}$. 
\begin{lemma} \label{lettim}  For $\alpha \in (0,1)$, there are constants $\varepsilon_0, c > 0$ such that
%\br \[
%\mathrm{P}(\tau^{x, \varepsilon}(S_{\varepsilon^\alpha}) > c_1  \varepsilon^{\alpha-1}) \geq 1/2, 
%\] \be \br ($<---$  maybe not needed) \be
\[
\mathrm{E} \tau^{x, \varepsilon}(S_{\varepsilon^\alpha}) \leq c  \varepsilon^{\alpha-1}
\]
for $x \in \overline{D}$, $\varepsilon \leq \varepsilon_0$.  For each $\delta > 0$ and $t(\varepsilon) \ll \varepsilon^{-1}$, 
\[
\mathrm{P}(\tau^{x, \varepsilon}(S_\delta) \leq t(\varepsilon)) \rightarrow 0
\]
as $\varepsilon \downarrow 0$. 
\end{lemma}
\proof All the statements easily follow from Lemmas~\ref{timeexit}, \ref{dirofexit}, and the strong Markov property of the process, once
we observe that the process coincides with the Brownian motion outside $\partial D$. \qed 
\\

Finally, we describe the location of the exit from a small neighborhood of $D$. While $X^{x,\varepsilon}_{\tau^{x, \varepsilon}(S_{\varepsilon^\alpha})}$ is
distributed on $S_{\varepsilon^\alpha}$, we can talk about the convergence of this distribution to a measure on $\partial D$, since $S_{\varepsilon^\alpha}$
can be viewed as a small perturbation of $\partial D$ as $\varepsilon \downarrow 0$. 
\begin{lemma} \label{stre}
For each $f \in C(\mathbb{T}^d)$, $\alpha \in (0,1)$,
\[
\lim_{\varepsilon \downarrow 0} \mathrm{E} f(X^{x,\varepsilon}_{\tau^{x, \varepsilon}(S_{\varepsilon^\alpha})}) = \int_{\partial D} f d \bar{\nu},
\]
uniformly in $x \in \overline{D}$, 
where $\bar{\nu}$ is the normalized Lebesgue measure on $\partial D$. 
\end{lemma} 
\proof For $x \in D \bigcup \Gamma^+_{\varepsilon^\alpha}$, consider the auxiliary process $\hat{X}^{x,\varepsilon}_t$ obtained from  ${X}^{x,\varepsilon}_t$ by
reflecting it (orthogonally to the surface) at $S_{\varepsilon^\alpha}$. As we have shown in Section~\ref{dttpA} for a similar process,
the measure $\mu^\varepsilon$, whose density with respect to the Lebesgue measure $\lambda$ is
\[
p^\varepsilon(x) = \left\{ \begin{array}{ll}
                 1, ~~~~~~ x \in \Gamma^+_{\varepsilon^\alpha},\\
                 \varepsilon^{-1}, ~~~ x \in D,
            \end{array}
            \right.
\]
is invariant for the family  $\hat{X}^{x,\varepsilon}_t$, $ x \in D \bigcup \Gamma^+_{\varepsilon^\alpha}$. 
Let the probability measure $\bar{\mu^\varepsilon}$
be obtained by multiplying $\mu^\varepsilon$ by a positive constant.

Let $\beta = (1 + \alpha)/2 \in (\alpha, 1)$. Take an arbitrary closed set 
$A \subset D$ with
a smooth boundary $\partial {A}$. We'll consider successive visits by the process  $\hat{X}^{x,\varepsilon}_t$ to $S_{\varepsilon^{\beta}}$ 
and $\partial {A}$. Namely, let $\tau^{x, \varepsilon}_0 = 0$,  
$\sigma^{x, \varepsilon}_n = \inf(t \geq \tau^{x, \varepsilon}_{n-1}: X^{x, \varepsilon}_t \in \partial A )$, $n \geq 1$,  while
$\tau^{x, \varepsilon}_n = \inf(t \geq \sigma^{x, \varepsilon}_n: X^{x, \varepsilon}_t \in    S_{\varepsilon^{\beta}}    )$, 
$n \geq 1$. 

Thus, for $x \in S_{\varepsilon^{\beta}}$, 
$\hat{X}^{x,\varepsilon}_{\tau^{x, \varepsilon}_n}$, $n \geq 0$,  is a Markov chain with the state space $S_{\varepsilon^{\beta}}$.
%For $x \in \partial A$, $\hat{X}^{x,\varepsilon}_{\sigma^{x, \varepsilon}_n}$, $n \geq 1$,  is a Markov chain with the state space $\partial A$. 
Let $\nu^\varepsilon$ 
%and $\eta^\varepsilon$ 
be the invariant measure for this chain.
%s for the two Markov chains, respectively. 
Let $\hat{f}^\varepsilon = f \chi_{\Gamma^+_{2\varepsilon^{\beta}} \setminus \Gamma^+_{\varepsilon^{\beta}}}$ (this function is non-zero in a thin
strip near $S_{\varepsilon^{\beta}}$). Then 
\begin{equation} \label{kha}
\int_{S_{\varepsilon^{\beta}}} \left(\mathrm{E} \int_0^{\tau^{x,\varepsilon}_1} \hat{f}^\varepsilon(\hat{X}^{x,\varepsilon}_t) dt \right) \nu^\varepsilon (dx) = 
\int_{S_{\varepsilon^{\beta}}} 
\mathrm{E} \tau^{x,\varepsilon}_1 \nu^\varepsilon (dx) \int_{\mathbb{T}^d} \hat{f}^\varepsilon d \bar{\mu}^\varepsilon.
\end{equation}  
Since $\hat{X}^{x,\varepsilon}_t$ coincides with the Brownian motion in the interior of $\Gamma^+_{\varepsilon^\alpha}$,
it is clear that there is $c_1(\varepsilon)$, independent of $f$, such that
\[
\lim_{\varepsilon \downarrow 0} \frac{\mathrm{E} \int_0^{\tau^{x,\varepsilon}_1} 
\hat{f}^\varepsilon(\hat{X}^{x,\varepsilon}_t) dt }{c_1(\varepsilon) f(x)} = 1,
\]
uniformly in $x \in S_{\varepsilon^{\beta}}$. Therefore, the left hand side of (\ref{kha}) is asymptotically equivalent to
\[
c_1(\varepsilon) \int_{S_{\varepsilon^{\beta}}} f d\nu^\varepsilon.
\]
The right hand side of (\ref{kha}) is asymptotically equivalent to 
$c_2(\varepsilon) \int_{\mathbb{T}^d} \hat{f}^\varepsilon d \bar{\mu}^\varepsilon$, with $c_2(\varepsilon)$ that is independent of $f$, which, in turn, is
asymptotically equivalent to $c_3(\varepsilon) \int_{\partial D} f d \bar{\nu}$, with $c_3(\varepsilon)$ that is independent of $f$. Thus
\[
\lim_{\varepsilon \downarrow 0}  \frac{c_1(\varepsilon)}{c_3(\varepsilon)} \int_{S_{\varepsilon^{\beta}}} f d\nu^\varepsilon  =
\int_{\partial D} f d \bar{\nu}.
\]
Since $c_1$ and $c_3$ do not depend on $f$,
\begin{equation} \label{nnii}
\lim_{\varepsilon \downarrow 0}  \int_{S_{\varepsilon^{\beta}}} f d \nu^\varepsilon  =
\int_{\partial D} f d \bar{\nu}.
\end{equation}
To complete the proof of the lemma, we write
\[
 \mathrm{E} f({X}^{x,\varepsilon}_{\tau^{x, \varepsilon}(S_{\varepsilon^\alpha})}) = 
 \mathrm{E} f(\hat{X}^{x,\varepsilon}_{\tau^{x, \varepsilon}(S_{\varepsilon^\alpha})}) = 
\]
\begin{equation} \label{mmd}
\mathrm{E} \left(f(\hat{X}^{x,\varepsilon}_{\tau^{x, \varepsilon}(S_{\varepsilon^\alpha})} )\chi_{\{\tau^{x, \varepsilon}(S_{\varepsilon^\alpha}) < 
\tau_1^{x, \varepsilon}\}}\right) + 
\end{equation}
\[
\sum_{n=1}^\infty \mathrm{E}
\left(f(\hat{X}^{x,\varepsilon}_{\tau^{x, \varepsilon}(S_{\varepsilon^\alpha})} )\chi_{\{ 
\tau_n^{x, \varepsilon} \leq \tau^{x, \varepsilon}(S_{\varepsilon^\alpha}) <
\tau_{n+1}^{x, \varepsilon}\}}\right). 
\]
The first term on the right hand side, as well as each individual term in the infinite sum, tends to zero when $\varepsilon \downarrow 0$, as 
follows from Lemma~\ref{dirofexit} and the strong Markov property of the process. In order to deal with the infinite sum, we write
\begin{equation} \label{lnln}
\mathrm{E}
\left(f(\hat{X}^{x,\varepsilon}_{\tau^{x, \varepsilon}(S_{\varepsilon^\alpha})} )\chi_{\{ 
\tau_n^{x, \varepsilon} \leq \tau^{x, \varepsilon}(S_{\varepsilon^\alpha}) <
\tau_{n+1}^{x, \varepsilon}\}}\right) =
\end{equation}
\[
\mathrm{E}
\left(f(\hat{X}^{x,\varepsilon}_{\tau^{x, \varepsilon}(S_{\varepsilon^\alpha})} ) \chi_{\{ 
\tau^{x, \varepsilon}(S_{\varepsilon^\alpha}) <
\tau_{n+1}^{x, \varepsilon}\}} | \tau_n^{x, \varepsilon} \leq \tau^{x, \varepsilon}(S_{\varepsilon^\alpha}) \right) \mathrm{P} 
\left( \tau_n^{x, \varepsilon} \leq \tau^{x, \varepsilon}(S_{\varepsilon^\alpha}) \right).
\]
Let the measure $\nu_n^{x,\varepsilon}$ on $S_{\varepsilon^\alpha}$ be defined, for Borel sets $A \subseteq S_{\varepsilon^\alpha}$, via
\[
\nu_n^{x,\varepsilon}(A) = \frac{\mathrm{P}(\hat{X}^{x,\varepsilon}_{\tau^{x, \varepsilon}_n} 
\in A,~\tau_n^{x, \varepsilon} \leq \tau^{x, \varepsilon}(S_{\varepsilon^\alpha}))}{\mathrm{P}(\tau_n^{x, \varepsilon} \leq \tau^{x, \varepsilon}(S_{\varepsilon^\alpha}))}.
\]
Then
\[
\mathrm{E}
\left(f(\hat{X}^{x,\varepsilon}_{\tau^{x, \varepsilon}(S_{\varepsilon^\alpha})} ) \chi_{\{ 
\tau^{x, \varepsilon}(S_{\varepsilon^\alpha}) <
\tau_{n+1}^{x, \varepsilon}\}} | \tau_n^{x, \varepsilon} \leq \tau^{x, \varepsilon}(S_{\varepsilon^\alpha}) \right) =
\]
\[
\int_{S_{\varepsilon^{\beta}}} \mathrm{E}
\left(f(\hat{X}^{y,\varepsilon}_{\tau^{y, \varepsilon}(S_{\varepsilon^\alpha})} ) \chi_{\{ 
\tau^{y, \varepsilon}(S_{\varepsilon^\alpha}) <
\tau_{1}^{y, \varepsilon}\}}  \right)  \nu_n^{x,\varepsilon}(d y).
\]
This quantity is asymptotically equivalent to $c_4(\varepsilon) \int_{S_{\varepsilon^{\beta}}} f d\nu^{x,\varepsilon}_n $ for some $c_4$ that does not
depend on $f$. 

Using the mixing properties of $\hat{X}^{x,\varepsilon}_t$, it is not difficult to show that $\nu_n^{x,\varepsilon}$ and $\nu^\varepsilon$ are close when $n$ is large and
$\varepsilon$ is small, in the sense that for each $\eta > 0$ there are $n_0$ and $\varepsilon_0$ such that
\[
|\int_{S_{\varepsilon^{\beta}}} f d\nu^{x,\varepsilon}_n - \int_{S_{\varepsilon^{\beta}}} f d\nu^\varepsilon| < \eta
\]
for all $n \geq n_0$, all $\varepsilon \leq \varepsilon_0$, and all $x \in \overline{D}$. Therefore, by (\ref{mmd}) and (\ref{lnln}),
$\mathrm{E} f({X}^{x,\varepsilon}_{\tau^{x, \varepsilon}(S_{\varepsilon^\alpha})})$ is asymptotically equivalent to
$c_5(\varepsilon) \int_{S_{\varepsilon^{\beta}}} f d\nu^\varepsilon$, where $c_5$ does not depend on $f$. In particular, applying this to $f \equiv 1$,
we obtain that $c_5 = 1$ works, i.e., 
\[
|\mathrm{E} f({X}^{x,\varepsilon}_{\tau^{x, \varepsilon}(S_{\varepsilon^\alpha})}) - \int_{S_{\varepsilon^{\beta}}} f d\nu^\varepsilon| \rightarrow 0.
\]
Combining this with (\ref{nnii}), we obtain the statement of the lemma.
\qed
\\
\\ 
{\bf Remark.} The assumption made in Lemmas~\ref{clex}-\ref{locexit} and in Lemma~\ref{stre}  that $D$ is a single trapping domain was notationally convenient, but not necessary 
for the results to hold (the proofs require only minor modifications). We can, therefore, use these lemmas in the general case.

\section{Behavior of the process inside of a trapping domain} \label{seconv}

%In this section we prove two lemmas: the first concerns the limiting behavior,  at finite time scales, inside of a trapping domain, while the second provides
%an estimate on the time it takes a process to leave a vicinity of a trapping domain. 

We focus on the case of a single trapping domain $D \subset \mathbb{T}^d$. 
Together with the family $X^{x, \varepsilon}_t$, we consider the family $Z^x_t$ constructed in Section~\ref{dttpA}. For $x \in \overline{D}$, $Z^x_t$ is just a Wiener process in $\overline{D}$
reflected at the boundary.
\begin{lemma} \label{lonref}
For each $x \in \mathbb{T}^d$, the processes $X^{x, \varepsilon}_t$ converge, as $\varepsilon \downarrow 0$, in distribution, to the process $Z^x_t$.
\end{lemma} 
\proof Let $A$ be the generator of $Z^x_t$. From the Hille-Yosida theorem it follows there is a dense linear subspace $\Psi$ of $C(\mathbb{T}^d)$ such that
for each  $\lambda > 0$  and  for  each $f \in \Psi$ the equation $\lambda F -  A F = f$ has a solution~$F \in D(A)$. By Lemma~\ref{ntch}, it is only remains to prove that for each $t > 0$ and each $f \in D(A)$,
  \[
	\lim_{\varepsilon \downarrow 0} \mathrm{E} (f( X^{x, \varepsilon}_t ) - f (  X^{x, \varepsilon}_0) - \int_0^t A f (  X^{x, \varepsilon}_u) du) =  0,
	\]
uniformly in $x \in \mathbb{T}^d$.

%Let $\Gamma^\varepsilon = \{x \in \overline{U}: {\rm dist}(x, \partial D) \leq \sqrt{\varepsilon} \}$ and $\gamma^\varepsilon = \{x \in \overline{U}: {\rm dist}(x, \partial D) = \sqrt{\varepsilon}\}$.
%For all sufficiently small $\varepsilon$, this is a smooth curve.
Let $\sigma^{x, \varepsilon}_0 = 0$,  $\tau^{x, \varepsilon}_n = \inf(t \geq \sigma^{x, \varepsilon}_{n-1}: X^{x, \varepsilon}_t \in \partial D)$, $n \geq 1$,  while
$\sigma^{x, \varepsilon}_n = \inf(t \geq \tau^{x, \varepsilon}_n: X^{x, \varepsilon}_t \in S_{-\sqrt{\varepsilon}} \bigcup S_{\sqrt{\varepsilon}})$, 
$n \geq 1$. Then 
\[
\mathrm{E} (f( X^{x, \varepsilon}_t ) - f (  X^{x, \varepsilon}_0) - \int_0^t A f (  X^{x, \varepsilon}_u) du) =
\]
\[
\mathrm{E} \sum_{n = 1}^\infty  \left(f(X^{x, \varepsilon}_{\tau^{x, \varepsilon}_n \wedge t} ) - f(X^{x, \varepsilon}_{\sigma^{x, \varepsilon}_{n-1} \wedge t}) - \frac{1}{2}  \int_{\sigma^{x, \varepsilon}_{n-1} \wedge t}^{\tau^{x, \varepsilon}_n \wedge t} \Delta f (X^{x, \varepsilon}_u) du \right) +
\]
\[
\mathrm{E} \sum_{n = 1}^\infty  \left(f(X^{x, \varepsilon}_{\sigma^{x, \varepsilon}_n \wedge t} ) - f(X^{x, \varepsilon}_{\tau^{x, \varepsilon}_{n} \wedge t}) - \frac{1}{2}  \int_{\tau^{x, \varepsilon}_{n} \wedge t}^{\sigma^{x, \varepsilon}_n \wedge t} \Delta f (X^{x, \varepsilon}_u) du \right).
\]
The first expectation on the right hand side is equal to zero since $X^{x, \varepsilon}_t$ is a Wiener process on $\mathbb{T}^d \setminus \partial D$. Our goal is to show that the second expectation tends to zero.
First, we need to control the number of terms in the sum.
Let us show that there is $c= c(t) > 0$ such that 
\begin{equation} \label{jj1}
 \mathrm{P}(\tau^{x, \varepsilon}_n  \leq t) \leq \exp({-{ c (n-1) \sqrt{\varepsilon}}}),~~~x \in {\mathbb{T}^d},~ n \geq 2.
\end{equation}
Since $\partial D$ is smooth, there is $r > 0$ such that the ball of radius $r$ tangent to $\partial D$ at $x$ lies either entirely in $\overline{U}$ or
entirely in $\overline{D}$. Let $\eta^\varepsilon$ be the time it takes a Wiener process starting inside a ball of radius $r$ at a distance $\sqrt{\varepsilon}$ from the boundary to reach the boundary. It is easy to see that there is $c = c(t) > 0$ such that
\begin{equation} \label{jj1s}
 \mathrm{P}(\eta^\varepsilon  \leq t) \leq \exp({-{ c  \sqrt{\varepsilon}}}).
\end{equation}
%Therefore, if $\eta^\varepsilon_k$, $k \geq 1$, is a sequence of such independent random variables, then
%\begin{equation} \label{jj1s}
% \mathrm{P}(\eta^\varepsilon_1 + ... + \eta^\varepsilon_n  \leq t) \leq \exp({-{c n \sqrt{\varepsilon}}}),~~~ n \geq 1.
%\end{equation}
By our construction, $\mathrm{P}(\tau^{x, \varepsilon}_n - \sigma^{x, \varepsilon}_{n-1} \leq z| X^{x, \varepsilon}_{\sigma^{x, \varepsilon}_{n-1}}) \leq
\mathrm{P}(\eta^\varepsilon \leq z)$ for each $n \geq 2$ and $z>0$. Therefore, (\ref{jj1}) follows from (\ref{jj1s}) and the strong Markov property. 
\\

Lemma \ref{timeexit}, together with (\ref{jj1}) and the strong Markov property of the process imply that
 there is $c= c(t)$ such that
\[
 \mathrm{E}  \sum_{n = 1}^\infty (\sigma^{x, \varepsilon}_{n} \wedge t - \tau^{x, \varepsilon}_n \wedge t) \leq c \sqrt{\varepsilon},~~~x \in {\mathbb{T}^d}.
\]
Therefore, since $\Delta f$ is bounded, 
\[
\lim_{\varepsilon \downarrow 0} \mathrm{E} 
\sum_{n = 1}^\infty \int_{\tau^{x, \varepsilon}_{n} \wedge t}^{\sigma^{x, \varepsilon}_n \wedge t} \Delta f (X^{x, \varepsilon}_u) du=0, 
\]
uniformly in $x \in \mathbb{T}^d$. It remains to show that 
\begin{equation} \label{rems}
 \lim_{\varepsilon \downarrow 0} \sup_{x \in \mathbb{T}^d} | \mathrm{E} \sum_{n = 1}^\infty  \left(f(X^{x, \varepsilon}_{\sigma^{x, \varepsilon}_n \wedge t} ) - f(X^{x, \varepsilon}_{\tau^{x, \varepsilon}_{n} \wedge t}) \right)| = 0.
\end{equation}
%Introduce the stopping time
%\[
%s'(t) = \Big \{
%  \begin{tabular}{ccc}
%  $\sigma^{x, \varepsilon}_n~~{\rm if}~~\tau^{x, \varepsilon}_{n} < t \leq \sigma^{x, \varepsilon}_n$ \\
%  $t~~~~~~~~~$ otherwise.
%  \end{tabular}
%\]
Now (\ref{rems}) will follow if we show that
\begin{equation} \label{rems2}
  \lim_{\varepsilon \downarrow 0} \sup_{x \in \mathbb{T}^d} | \mathrm{E} \sum_{n = 1}^\infty  \chi_{\{\tau^{x,\varepsilon}_n < t\}} 
	\left(f(X^{x, \varepsilon}_{\sigma^{x, \varepsilon}_n } ) - f(X^{x, \varepsilon}_{\tau^{x, \varepsilon}_{n} }) \right) | = 0,
\end{equation}
since the difference between (\ref{rems2}) and (\ref{rems}) is estimated from above by $\sup_{x,y \in \Gamma_{\sqrt{\varepsilon}}}|f(x) - f(y)|$,
 which goes to zero as $\varepsilon \downarrow 0$. Let $N^{x,\varepsilon} = \max(n: \tau^{x, \varepsilon}_n < t) $.  By the strong Markov property,
\[
\sup_{x \in \mathbb{T}^d} | \mathrm{E} \sum_{n = 1}^\infty  \chi_{\{\tau^{x,\varepsilon}_n < t\}} 
	\left(f(X^{x, \varepsilon}_{\sigma^{x, \varepsilon}_n } ) - f(X^{x, \varepsilon}_{\tau^{x, \varepsilon}_{n} }) \right) | \leq
\sup_{x \in \mathbb{T}^d} \mathrm{E} N^{x,\varepsilon}
\sup_{x \in \partial D} |\mathrm{E} (f(X^{x, \varepsilon}_{\sigma^{x, \varepsilon}_1}) - f(x)) |.
\]
Since $\langle \nabla f_{\overline{D}} (x) , n_D(x) \rangle = 0$ for $x \in \partial D$, 
the right hand side tends to zero by (\ref{jj1}) and Lemma~\ref{locexit}. This concludes the proof of Lemma~\ref{lonref}.
\qed
\\
\\
{\bf Remark:}  Before Lemma~\ref{lonref}, we made the assumption that $D$ is the only trapping domain. In fact, it is clear that the result holds even if
there are other trapping domains, as long as they are disjoint from $D$, and the initial point $x$ belongs to $\overline{D}$. 

%{\bf Remark on notation.} From this point on, we'll use the notation ${\alpha}^{\varepsilon}(\delta)$ for a generic quantity that satisfies  
%\[
%\lim_{\delta \downarrow 0} \limsup_{\varepsilon \downarrow 0}  |{\alpha}^{\varepsilon}(\delta)| = 0.
%\] 
%The notation may stand for different quantities from line to line. 

%Let $D^\delta = \{ x \in \mathbb{T}^d, {\rm dist}(x, D) < \delta \}$ be the $\delta$-neighborhood of $D$. For a closed set $A \subseteq \mathbb{T}^d$, let
%\[
%\tau^{x, \varepsilon}(A) = \inf\{t \geq 0: X^{x, \varepsilon}_t \in A\}. 
%\]
%\begin{lemma}
%There is a constant $c_D > 0$ such that
%\[
%\tau^{x, \varepsilon}(\partial D^\delta)  = \varepsilon \delta (c_D + {\alpha}^{x, \varepsilon}(\delta)),
%\]
%where
%\[
%\lim_{\delta \downarrow 0} \limsup_{\varepsilon \downarrow 0}  |{\alpha}^{x,\varepsilon}(\delta)| = 0,
%\] 
% uniformly for  $x \in \overline{D}$.
%\end{lemma}

 \begin{lemma} \label{cor82} Suppose that $1 \ll t_1(\varepsilon) \leq t_2(\varepsilon) \ll \varepsilon^{-1}$. Then
\[
\inf_{x \in \overline{D}, t_1(\varepsilon) \leq t \leq t_2(\varepsilon)}  \mathrm{P}(X^{x, \varepsilon}_{t} \in \overline{D}) \rightarrow 1
\]
as $\varepsilon \downarrow 0$. 
\end{lemma}
 \proof   It is not difficult to show that the convergence in Lemma~\ref{lonref} is uniform with respect to the initial point. In particular,
for each $f \in C(\mathbb{T}^d)$,  $\mathrm{E} f(X^{x, \varepsilon}_t) \rightarrow \mathrm{E} f ( Z^x_t)$  
as $\varepsilon \downarrow 0$ uniformly in $x \in \mathbb{T}^d$.   For $ \eta > 0$, choose $\delta > 0$ and 
$0 \leq f \leq 1$ such that ${\rm supp}(f) \subset D$ and 
$ \mathrm{E} f ( Z^x_1) \geq 1 - \eta$ whenever ${\rm dist} (x, D) \leq \delta$. From the second statement in Lemma~\ref{lettim} and the Markov property
applied to time $t-1$, if follows that $ \mathrm{P}(X^{x, \varepsilon}_{t} \in \overline{D}) \geq 1 - 2\eta$ for all sufficiently small $\varepsilon$,
which gives the desired result. 
\qed 
\\

\section{The limiting behavior of the  trace process} \label{seconv2}

Let us assume that
$S$ contains all the indices $l$ such that $D_l \prec D_0$. Let $Y^{x, \varepsilon}_t $, $x \in \mathbb{T}^d$, be the $U'$-valued processes defined
in Section~\ref{abstle} and $Y^x_t$, $x \in U'$, be the  $U'$-valued processes defined
in Section~\ref{dttp} corresponding to $k = 0$ and $S$.
\begin{lemma} \label{mt1xx}
For each $x \in  \mathbb{T}^d$, the measures on $C([0,\infty), U')$ induced by the processes $Y^{x, \varepsilon}_t $ converge weakly, as $\varepsilon \downarrow 0$, to the measure induced by $Y^{h(x)}_t$. 
\end{lemma}
\proof Let $L$ be the generator of $Y^x_t$. 
  As in the proof of Lemma~\ref{lonref} (now using Lemma~\ref{flei} instead of Lemma~\ref{ntch}), it is 
sufficient to show that for each $t > 0$ and each $f \in D(L)$,
\[
\lim_{\varepsilon \downarrow 0} \mathrm{E} \left(f(Y^{x, \varepsilon}_t) - f(x) -  \frac{1}{2} \int_0^t \Delta f (Y^{x, \varepsilon}_u) du \right) = 0,
\]
 uniformly in  $x \in \overline{U}$.

Again, we define two sequences of stopping times, but somewhat differently from the way it was done in Section~\ref{seconv}.
Let $\sigma^{x, \varepsilon}_0 = 0$,  $\tau^{x, \varepsilon}_n = \inf(t \geq \sigma^{x, \varepsilon}_{n-1}: X^{x, \varepsilon}_t \in \partial D)$, $n \geq 1$, 
while $\sigma^{x, \varepsilon}_n = \inf(t \geq \tau^{x, \varepsilon}_n: X^{x, \varepsilon}_t \in S_{\sqrt{\varepsilon}})$, $n \geq 1$. Then 
\[
\mathrm{E} \left(f(Y^{x, \varepsilon}_t) - f(x) -  \frac{1}{2}  \int_0^t \Delta f (Y^{x, \varepsilon}_u) du \right) =
\]
\[
\mathrm{E} \sum_{n = 1}^\infty  \left(f(X^{x, \varepsilon}_{\tau^{x, \varepsilon}_n \wedge s(t)} ) - f(X^{x, \varepsilon}_{\sigma^{x, \varepsilon}_{n-1} \wedge s(t)}) - \frac{1}{2}  \int_{\sigma^{x, \varepsilon}_{n-1} \wedge s(t)}^{\tau^{x, \varepsilon}_n \wedge s(t)} \Delta f (X^{x, \varepsilon}_u) du \right) +
\]
\[
\mathrm{E} \sum_{n = 1}^\infty  \left(f(X^{x, \varepsilon}_{\sigma^{x, \varepsilon}_n \wedge s(t)} ) - f(X^{x, \varepsilon}_{\tau^{x, \varepsilon}_{n} \wedge s(t)}) - \frac{1}{2}  \int_{\tau^{x, \varepsilon}_{n} \wedge s(t)}^{\sigma^{x, \varepsilon}_n \wedge s(t)} \Delta f (X^{x, \varepsilon}_u) du \right),
\]
where we put $\Delta f \equiv 0 $ on $D$. The first expectation on the right hand side is equal to zero since $X^{x, \varepsilon}_t$ is a Wiener process on $\overline{U}$. Our goal is to show that the second expectation tends to zero. 
As in the proof of Lemma~\ref{lonref} (but now with $s(t)$ instead of $t$), 
there is $c= c(t) > 0$ such that 
\begin{equation} \label{jj1x}
 \mathrm{P}(\tau^{x, \varepsilon}_n  \leq s(t)) \leq \exp({-{ c (n-1) \sqrt{\varepsilon}}}),~~~x \in {\mathbb{T}^d},~ n \geq 2.
\end{equation}

Lemma \ref{stwo}, together with (\ref{jj1x}) and the strong Markov property of the process imply that
 there is $c= c(t)$ such that
\[
 \mathrm{E}  \sum_{n = 1}^\infty \lambda(\{u:  \tau^{x, \varepsilon}_{n} \wedge s(t) \leq u \leq \sigma^{x, \varepsilon}_n \wedge s(t),~~X^{x, \varepsilon}_u \in 
\Gamma^+_{\sqrt{\varepsilon}} \}) \leq c \sqrt{\varepsilon},~~~x \in {\mathbb{T}^d}.
\]
Therefore, since $\Delta f$ is bounded, 
\[
\lim_{\varepsilon \downarrow 0} \mathrm{E} 
\sum_{n = 1}^\infty \int_{\tau^{x, \varepsilon}_{n} \wedge s(t)}^{\sigma^{x, \varepsilon}_n \wedge s(t)} \Delta f (X^{x, \varepsilon}_u) du=0, 
\]
uniformly in $x \in \mathbb{T}^d$. It remains to show that
\begin{equation} \label{remsxx}
 \lim_{\varepsilon \downarrow 0} \sup_{x \in \mathbb{T}^d} | \mathrm{E} \sum_{n = 1}^\infty  \left(f(X^{x, \varepsilon}_{\sigma^{x, \varepsilon}_n \wedge s(t)} ) - f(X^{x, \varepsilon}_{\tau^{x, \varepsilon}_{n} \wedge s(t)}) \right)| = 0.
\end{equation}
%Introduce the stopping time
%\[
%s'(t) = \Big \{
%  \begin{tabular}{ccc}
%  $\sigma^{x, \varepsilon}_n~~{\rm if}~~\tau^{x, \varepsilon}_{n} < t \leq \sigma^{x, \varepsilon}_n$ \\
%  $t~~~~~~~~~$ otherwise.
%  \end{tabular}
%\]
As in the proof of Lemma~\ref{lonref}, (\ref{remsxx}) will follow if we show that
\begin{equation} \label{rems2xx}
  \lim_{\varepsilon \downarrow 0} \sup_{x \in \mathbb{T}^d}  |\mathrm{E} \sum_{n = 1}^\infty  \chi_{\{\tau^{x,\varepsilon}_n < s(t)\}} 
	\left(f(X^{x, \varepsilon}_{\sigma^{x, \varepsilon}_n } ) - f(X^{x, \varepsilon}_{\tau^{x, \varepsilon}_{n} }) \right) | = 0,
\end{equation}
since the difference between (\ref{rems2xx}) and (\ref{remsxx}) is estimated from above by $\sup_{x,y \in \Gamma^+_{\sqrt{\varepsilon}}}|f(x) - f(y)|$,
 which goes to zero as $\varepsilon \downarrow 0$. Let $N^{x,\varepsilon} = \max(n: \tau^{x, \varepsilon}_n < s(t)) $.  By the strong Markov property,
\[
\sup_{x \in \mathbb{T}^d} | \mathrm{E} \sum_{n = 1}^\infty  \chi_{\{\tau^{x,\varepsilon}_n < s(t)\}} 
	\left(f(X^{x, \varepsilon}_{\sigma^{x, \varepsilon}_n } ) - f(X^{x, \varepsilon}_{\tau^{x, \varepsilon}_{n} }) \right) | \leq
\sup_{x \in \mathbb{T}^d} \mathrm{E} N^{x,\varepsilon}
\sup_{x \in \partial D} |\mathrm{E} (f(X^{x, \varepsilon}_{\sigma^{x, \varepsilon}_1}) - f(x)) |.
\]
The right hand side tends to zero by (\ref{jj1x}) and Lemma~\ref{stre}. This concludes the proof of Lemma~\ref{mt1xx}.
\qed
\\

Both Lemma~\ref{lonref} and Lemma~\ref{mt1xx} are somewhat restrictive: in the former, it is assumed that there are no sub-domains inside the trapping domain
under consideration, while
in the latter, the limiting process is considered in the space $U'$ corresponding to $k = 0$ rather than general $k$. These two lemmas can be combined, however, to treat the general case. 

Namely, let us assume that
$S$ contains all the indices $l$ such that $D_l \prec D_k$.
If $k \neq 0$, we need to distinguish between the spaces  $U_0 = D_0 \setminus \bigcup_{l \in S} \overline{D}_l$ and 
$U = D_k \setminus \bigcup_{l \in S} \overline{D}_l$.  As in Section~\ref{abstle}, we can define the mapping $h: \mathbb{T}^d \rightarrow U'_0$ 
by $h(x) = x$ for $x \in U_0$ and $h(x) = d_l$ for $x \in \overline{D}_l$. We can also define
process $Y^{x, \varepsilon}_t $, $x \in \mathbb{T}^d$, 
which is now $U'_0$-valued and is obtained from $X^{x, \varepsilon}_t $ by running the clock only when $X^{x, \varepsilon}_t $ is in $\overline{U}_0$.
On the other hand, the process  $Y^x_t$, $x \in U'$,  corresponding to $k$ and $S$ and defined in Section~\ref{dttp},  is $U'$-valued. Since $U' \subseteq
U'_0$, we can also view $Y^x_t$ as $U'_0$-valued.
\begin{lemma} \label{mt1xx22}
For each $x \in  \overline{D}_k$, the measures on $C([0,\infty), U'_0)$ induced by the processes $Y^{x, \varepsilon}_t $ converge weakly, as 
$\varepsilon \downarrow 0$, to the measure induced by $Y^{h(x)}_t$. 
\end{lemma}
\proof Assume that $k \neq 0$ (otherwise, the statement follows from Lemma~\ref{mt1xx}). Let $R$ be a domain with a smooth boundary such that
$\overline{D}_l \subset R$ for $l \in S$ and $\overline{R} \subset D_k$. Assume that $x \in \overline{D}_k \setminus R$ (the case when $x \in R$ is
treated similarly). Define two sequences of stopping times: $\sigma^{x, \varepsilon}_0 = 0$,  $\tau^{x, \varepsilon}_n = 
\inf(t \geq \sigma^{x, \varepsilon}_{n-1}: Y^{x, \varepsilon}_t \in \partial R)$, $n \geq 1$, while  $\sigma^{x, \varepsilon}_n = \inf(t \geq \tau^{x, \varepsilon}_n: Y^{x, \varepsilon}_t \in \partial D_k)$, $n \geq 1$. Also define $\overline{\sigma}^x_n$, $n \geq 0$, and $\overline{\tau}^x_n$, $n \geq 1$, in the same way,
but with $Y^{h(x)}_t$ instead of $Y^{x, \varepsilon}_t$.
From Lemma~\ref{lonref}, it follows that the measure induced by $Y^{x, \varepsilon}_{\tau^{x, \varepsilon}_1 \wedge t}$ 
converges weakly, as  $\varepsilon \downarrow 0$, to the measure induced by $Y^{h(x)}_{\overline{\tau}^x_1 \wedge t}$. From Lemma~\ref{mt1xx}, by the strong Markov property of the processes, it follows that the measure induced by $Y^{x, \varepsilon}_{\sigma^{x, \varepsilon}_1 \wedge t}$ 
converges to the measure induced by $Y^{h(x)}_{\overline{\sigma}^x_1 \wedge t}$. Continuing by induction, we obtain that,
for each $n$, the measure induced by $Y^{x, \varepsilon}_{\sigma^{x, \varepsilon}_n \wedge t}$ 
converges to the measure induced by $Y^{h(x)}_{\overline{\sigma}^x_n \wedge t}$. This implies the statement of the lemma.
\qed
\\

Let $f$ be a continuous function on $U'_0$ and let $t > 0$. From Lemma~\ref{mt1xx22}, it follows that for each	$x \in \overline{D}_k$, $\mathrm{E} 
f(Y^{x, \varepsilon}_t) \rightarrow  \mathrm{E} f(Y^{h(x)}_t)$ as $\varepsilon \downarrow 0$. 
Let $Q$ be a domain with smooth boundary such that $\overline{Q}  \subset D_k \setminus \bigcup_{l \in S} \overline{D}_l$. Let $\tau^{x, \varepsilon}$ be
the first time when $Y^{x, \varepsilon}_t$ reaches $\overline{Q}$ and $\tau^x$ be the first time when $Y^{h(x)}_t$  reaches $\overline{Q}$. Let $g \in C(\overline{Q})$.
From Lemma~\ref{mt1xx22}, it follows that for each	$x \in \overline{D}_k$, $\mathrm{E} 
g(Y^{x, \varepsilon}_{\tau^{x,\varepsilon}}) \rightarrow  \mathrm{E} g(Y^{h(x)}_{\tau^x})$ as $\varepsilon \downarrow 0$.
	It is not difficult to show that $\mathrm{E} 
f(Y^{x, \varepsilon}_t)$ and $\mathrm{E} 
g(Y^{x, \varepsilon}_{\tau^{x,\varepsilon}})$ are uniformly continuous in $x \in \overline{D}_k$ and $\varepsilon \leq \varepsilon_0$ for some $\varepsilon_0$.  Therefore, we have the following corollary.
\begin{corollary} \label{coco}
For each $f \in C(U'_0)$ and $t >0$, 
\[
\lim_{\varepsilon \downarrow 0}  \mathrm{E} f(Y^{x, \varepsilon}_t) =  \mathrm{E} f(Y^{h(x)}_t),
\]
uniformly in $x \in \overline{D}_k$. For each $g \in C(\overline{Q})$, 
\[
\lim_{\varepsilon \downarrow 0} 
\mathrm{E}  g(Y^{x, \varepsilon}_{\tau^{x,\varepsilon}}) =  \mathrm{E} g(Y^{h(x)}_{\tau^x}),
\]
uniformly in $x \in \overline{D}_k$.
\end{corollary}

\section{Proof of the main result} \label{pmree}
In this section, we prove Theorem~\ref{mteor} in the particular case outlined in the Introduction (see Figure~\ref{secondkind}).
We assume that $\varepsilon_k(\varepsilon) =
\varepsilon$ for each $1 \leq k \leq 7$. In this case, Assumption 1 of Section~\ref{mremre} holds. We distinguish different cases for the behavior
of $t(\varepsilon)$, each of which conforms with Assumption 2. 

First, consider the case when $1 \ll t(\varepsilon) \ll \varepsilon^{-1}$ and the
process starts in $\overline{D}_1$. To stress that the limits below are uniform in the choice of $t(\varepsilon)$, we consider
$1 \ll t_1(\varepsilon) \leq t_2(\varepsilon) \ll \varepsilon^{-1}$ and assume that $t_1(\varepsilon) \leq  t(\varepsilon) \leq t_2(\varepsilon)$.

For $x \in \overline{D}_1$, let 
$Y^{x}_t$ be the process corresponding to $k = 1$ and $S = \O $ in the notation of Section~\ref{dttp} (it is simply
the Brownian motion with instantaneous reflection on $\partial D_1$ to 
the interior of $D_1$). Let $f \in C(\mathbb{T}^d)$ and $\eta > 0$.
%Take $\delta(\varepsilon)$ such that $\varepsilon t(\varepsilon) \ll \delta(\varepsilon) \ll 1$. 
Take $s$ sufficiently large so that 
\begin{equation} \label{juju1}
|\mathrm{E} f(Y^x_{s}) - \int_{D_1} f d\lambda_1| \leq \eta/3
\end{equation}
for $x \in \overline{D}_1$,
% satisfying ${\rm dist}(x, D_1) \leq \delta(\varepsilon)$ and all sufficiently small $\varepsilon$ 
where $\lambda_1$ is the normalized Lebesgue measure on $D_1$. 
From  Lemma~\ref{cor82}, it follows that 
\begin{equation} \label{juju2}
\mathrm{P}(X^{x,\varepsilon}_{t(\varepsilon) - s} \in  \overline{D}_1) \geq 1 - \eta/3
\end{equation}
for all sufficiently small $\varepsilon$ and $x \in \overline{D}_1$. Combining (\ref{juju1}), (\ref{juju2}), using the Markov property of 
the process $X^{x,\varepsilon}_t$ (with the time $t(\varepsilon) - s$), and Corollary~\ref{coco}, we obtain that
\begin{equation} \label{clll}
|\mathrm{E} f(X^{x,\varepsilon}_{t(\varepsilon)}) - \int_{D_1} f d\lambda_1| \leq \eta,~~~x \in \overline{D}_1,~~~t_1(\varepsilon) \leq t(\varepsilon) \leq t_2(\varepsilon).
\end{equation}
This gives the desired result for $x \in \overline{D}_1$.
 Similarly, if
$x \in \overline{D}_k$ with $k = 2,3$, then the distribution of $X^{x, \varepsilon}_{t(\varepsilon)}$ will be asymptotically close to  $\lambda_k$.

Now consider the case when the process starts in $\overline{D}_4 \setminus \overline{D}_1$. 
 Let
$\tau^{x, \varepsilon} = \tau^{x, \varepsilon}(\partial D_1)$ be the first time when the process  $X^{x, \varepsilon}_{t}$
hits $\partial D_1$.  Since $X^{x}_t$ coincides with the Brownian motion in ${D}_4 \setminus \overline{D}_1$, from Lemma~\ref{dirofexit}
it follows that
\begin{equation} \label{shrr}
\mathrm{P}(\tau^{x, \varepsilon} \leq \frac{t(\varepsilon)}{2}) \rightarrow 1~~~{\rm as}~~\varepsilon \downarrow 0.
\end{equation}
Therefore, using (\ref{shrr}) and the strong Markov property of $X^{x, \varepsilon}_{t}$ (with the stopping time $\tau^{x, \varepsilon}$), we can conclude that (\ref{clll}) holds for $x \in \overline{D}_4 \setminus \overline{D}_1$.
The same argument, but with two stopping times, $X^{x, \varepsilon}_{t}$ hitting $\partial D_4$ and then hitting $\partial D_1$, lead to (\ref{clll}) for $x \in \overline{D}_6 \setminus \overline{D}_4$. Thus
\begin{equation} \label{clll2}
|\mathrm{E} f(X^{x,\varepsilon}_{t(\varepsilon)}) - \int_{D_1} f d\lambda_1| \leq \eta,~~~x \in \overline{D}_6,~~~t_1(\varepsilon) \leq t(\varepsilon) \leq t_2(\varepsilon).
\end{equation}
Similarly,
\begin{equation} \label{clll3}
|\mathrm{E} f(X^{x,\varepsilon}_{t(\varepsilon)}) - \int_{D_2} f d\lambda_2| \leq \eta,~~~x \in \overline{D}_5,~~~t_1(\varepsilon) \leq t(\varepsilon) \leq t_2(\varepsilon),
\end{equation}
and, as we already saw,
\begin{equation} \label{clll4}
|\mathrm{E} f(X^{x,\varepsilon}_{t(\varepsilon)}) - \int_{D_3} f d\lambda_3| \leq \eta,~~~x \in \overline{D}_3,~~~t_1(\varepsilon) \leq t(\varepsilon) \leq t_2(\varepsilon).
\end{equation}

Next, consider the case when $x \in \overline{D}_7 \setminus (\overline{D}_6 \bigcup \overline{D}_5 \bigcup \overline{D}_3)$. 
Consider the Wiener process with 
reflection on $\partial D_7$ ($Y^{x}_t$ corresponding to $k = 7$ and $S = \O $). Let
$\tau^{x, \varepsilon} = \tau^{x, \varepsilon}(\partial D_6 \bigcup \partial D_5 \bigcup \partial D_3)$ be the first time when the process  $X^{x, \varepsilon}_{t}$
hits $\partial D_6 \bigcup \partial D_5 \bigcup \partial D_3$. We define $\tau^x$ similarly, but for the process $Y^x_t$ rather than $X^{x,\varepsilon}_t$. 

Observe that (\ref{shrr}) holds with this new stopping time $\tau^{x, \varepsilon}$
and, by Corollary~\ref{coco}, 
\[
\mathrm{P}(X^{x,\varepsilon}_{\tau^{x,\varepsilon}} \in \partial D_6) \rightarrow \mathrm{P}(Y^{x}_{\tau^{x}} \in \partial D_6),
\]
\[
\mathrm{P}(X^{x,\varepsilon}_{\tau^{x,\varepsilon}} \in \partial D_5) \rightarrow \mathrm{P}(Y^{x}_{\tau^{x}} \in \partial D_5),
\]
\[
\mathrm{P}(X^{x,\varepsilon}_{\tau^{x,\varepsilon}} \in \partial D_3) \rightarrow \mathrm{P}(Y^{x}_{\tau^{x}} \in \partial D_3).
\]
Using the strong Markov property of $X^{x,\varepsilon}_t$ (with the stopping time $\tau^{x,\varepsilon}$), from (\ref{clll2})-(\ref{clll4}) we conclude that
 the distribution of $X^{x, \varepsilon}_{t(\varepsilon)}$ tends to 
\[
\mathrm{P}(Y^{x}_{\tau^{x}} \in \partial D_6) \lambda_1 + \mathrm{P}(Y^{x}_{\tau^{x}} \in \partial D_5) \lambda_2 +\mathrm{P}(Y^{x}_{\tau^{x}} \in \partial D_3)\lambda_3.
\]

If $x \in D_0 \setminus \overline{D}_7$, we consider the Wiener process in $D_0$ until the first time it hits $\partial D_7$. 
Let $\pi(x, y)$, $x \in D_0 \setminus \overline{D}_7$, $y \in \partial D_7$, be the corresponding Poisson kernel. We apply the arguments above, but starting with the
point $y$ where the process first reaches $\partial D_7$. Thus the distribution of $X^{x, \varepsilon}_{t(\varepsilon)}$ tends to 
\begin{equation} \label{innip}
\int_{\partial D_7} (\mathrm{P}(Y^{y}_{\tau^{y}} \in \partial D_6) \lambda_1 + \mathrm{P}(Y^{y}_{\tau^{y}} \in \partial D_5) \lambda_2 +\mathrm{P}(Y^{y}_{\tau^{y}} \in \partial D_3)\lambda_3)  \pi(x, y) dy.
\end{equation}

Now let us consider the case when $\varepsilon^{-1} \ll t(\varepsilon) \ll \varepsilon^{-2}$. Assume that $x \in \overline{D}_6$.   The transition probabilities between $\overline{D}_6$, $\overline{D_4}$, and
$\overline{D}_1$ are controlled by Lemma~\ref{dirofexit} and the transition times are at least of order one, since $X^{x, \varepsilon}_{t}$
coincides with the Brownian motion away from the boundaries of the domains. From here 
it easily follows that $X^{x, \varepsilon}_{t}$ does not leave the $\sqrt{\varepsilon}$-neighborhood of $\overline{D}_6$ in prior to 
time $t(\varepsilon)$ with probability that tends to one.
 Take an arbitrary $\tilde{t}(\varepsilon)$ such that $1 \ll \tilde{t}(\varepsilon)
\ll \varepsilon^{-1}$.  Let $\tau^{x, \varepsilon}$ be the first time after ${t}(\varepsilon) - \tilde{t}(\varepsilon)$ when $X^{x, \varepsilon}_{t} \in
\overline{D}_6$. Then
$\mathrm{P}(\tau^{x, \varepsilon} \leq {t}(\varepsilon) - \frac{\tilde{t}(\varepsilon)}{2}) \rightarrow 1$ as $\varepsilon \downarrow 0$. Using the strong Markov
property of the process $X^{x, \varepsilon}_{t}$ (with the stopping time $\tau^{x, \varepsilon}$) and the result on the limiting distribution of the process
at time scales that satisfy $1 \ll t(\varepsilon) \ll \varepsilon^{-1}$, we obtain that the distribution of $X^{x, \varepsilon}_{t(\varepsilon)}$ tends to $\lambda_1$.
Similarly, for $x \in \overline{D}_5$, the distribution of $X^{x, \varepsilon}_{t(\varepsilon)}$ tends to $\lambda_2$. 

For $x \in \overline{D}_7 \setminus (\overline{D_6} \cup \overline{D}_5)$ (including $x \in \overline{D}_3$), 
consider the process $Y^{h(x)}_t$  corresponding to $k = 7$ and $S = \{3\} $.
The process  $X^{x,\varepsilon}_t$ reaches 
$\partial{D}_6 \cup \partial{D}_5$ prior to time $t(\varepsilon)/2$ with probability that tends to one, as follows from the first statement of Lemma~\ref{lettim} and  Corollary~\ref{coco}.  Thus the limiting distribution of $X^{x, \varepsilon}_{t(\varepsilon)}$ is the weighted sum of $\lambda_1$ and
$\lambda_2$. From Corollary~\ref{coco}, it follows that, for $x \in \overline{D}_7 \setminus (\overline{D}_6 \cup \overline{D}_5)$,
\[
\lim_{\varepsilon \downarrow 0} \mathrm{P}(X^{x,\varepsilon}_{\tau^{x, \varepsilon}} \in \partial D_6) = 
\mathrm{P}(Y^{x}_{\tau^{x}} \in \partial D_6),~~\lim_{\varepsilon \downarrow 0} \mathrm{P}(X^{x,\varepsilon}_{\tau^{x, \varepsilon}} \in \partial D_5) = 
\mathrm{P}(Y^{x}_{\tau^{x}} \in \partial D_5),
\]
where
$\tau^{x, \varepsilon} = \tau^{x, \varepsilon}(\partial D_6 \bigcup \partial D_5)$ is the first time when $X^{x,\varepsilon}_t$
hits $ \partial D_6 \bigcup \partial D_6$ and the stopping time $\tau^x$ is defined similarly, with $Y^{h(x)}_t$ instead of $X^{x,\varepsilon}_t$.
Using the strong Markov
property of the process $X^{x, \varepsilon}_{t}$ (with the stopping time $\tau^{x, \varepsilon}$), we conclude that the 
distribution of $X^{x, \varepsilon}_{t(\varepsilon)}$ tends to 
\[
\mathrm{P}(Y^{x}_{\tau^{x}} \in \partial D_6) \lambda_6 + \mathrm{P}(Y^{x}_{\tau^{x}} \in \partial D_5) \lambda_5.
\]
If $x \in D_0 \setminus \overline{D}_7$, we can use the same arguments that led to (\ref{innip}) and obtain that
the distribution of $X^{x, \varepsilon}_{t(\varepsilon)}$ tends to 
\[
\int_{\partial D_7} (\mathrm{P}(Y^{y}_{\tau^{y}} \in \partial D_6) \lambda_1 + \mathrm{P}(Y^{y}_{\tau^{y}} \in \partial D_5) \lambda_2)  \pi(x, y) dy.
\]
Here, $\pi(x, y)$, $x \in D_0 \setminus \overline{D}_7$, $y \in \partial D_7$, is the same Poisson kernel as above, but  $Y^{x}_t$ is 
the process corresponding to $k = 7$ and $S = \{3\} $.

The case when $\varepsilon^{-2} \ll t(\varepsilon) \ll \varepsilon^{-3}$,  $X^{x, \varepsilon}_{t}$ is similar. The process $X^{x,\varepsilon}_t$
has enough time to reach $\overline{D}_6$ from each $x$, but not enough time to exit a $\sqrt{\varepsilon}$ neighborhood of a $D_6$. Thus, for all $x$,  $X^{x, \varepsilon}_{t(\varepsilon)}$ tends to 
a uniform distribution on $D_1$. 

In the cases when $\varepsilon^{-3} \ll t(\varepsilon) \ll \varepsilon^{-4}$, and $t(\varepsilon) \gg \varepsilon^{-4}$, we can take an arbitrary function
$\tilde{t}$ such that $\varepsilon^{-2} \ll \tilde{t}(\varepsilon) \ll \varepsilon^{-3}$. Using the Markov property of $X^{x,\varepsilon}_t$ (with time $t(\varepsilon)-
\tilde{t}(\varepsilon)$) and the fact that $X^{x, \varepsilon}_{\tilde{t}(\varepsilon)}$ tends to 
a uniform distribution on $D_1$, we obtain that $X^{x, \varepsilon}_{t(\varepsilon)}$ tends to 
a uniform distribution on $D_1$. 
\\
\\\\
\noindent {\bf \large Acknowledgments}:  While working on this
article, L. Koralov was supported by the ARO grant W911NF1710419.
\\
\\


\begin{thebibliography}{999999}

%\bibitem{AF} Athreya A., Freidlin M.I., {\it Metastability and
%Stochastic Resonance in Nearly-Hamiltonian Systems}, Stochastics
%and Dynamics, 8, 1, pp 1-21, 2008.

%\bibitem{CF} Chen Z., Freidlin M.I., {\it Smoluchowski-Kramers
%Approximation and Exit Problems}, Stochastics and Dynamics, 5, pp
%569-585, 2005.
%
%\bibitem{DSt} Dueschel J-D., Strook D., {Large Deviations}, AMS
%Chelsea Publishing, 2001.

%\bibitem{F} Freidlin M. I., {\it Metastability and Stochastic
%Resonance of Multiscale Systems}, Contemporary Mathematics, 2008.

%\bibitem{F77} Freidlin M. I., {\it Sublimiting Distributions and
%Stabilization of Solutions of Parabolic Equations with a Small
%Parameter}, Soviet Math. Dokl., 235, 5, pp 1042-1045, 1977.

%\bibitem{F85} Freidlin M.I., {\it Functional Integration and
%Partial Differential Equations}, Princeton University Press, 1985.

%\bibitem{FPhD} Freidlin M. I., {\it Quasi-deterministic
%Approximation, Metastability and Stochastic Resonance}, Physica D,
%137, pp 333-352, 2000.

%\bibitem{FKo} Freidlin M. I., Koralov L. {\it Metastability for Nonlinear Random Perturbations of Dynamical Systems}, Stochastic
%Processes and Applications 120 (2010), no. 7, 1194–-1214.

%\bibitem{FK} Freidlin M. I., Koralov L. {\it Nonlinear Stochastic Perturbations of Dynamical Systems and
%Quasi-linear Parabolic PDE's with a Small Parameter}, Probability
%Theory and Related Fields (2010), 147, pp 273-301. (An updated
%version is available on the arXiv.)

%\bibitem{FWa} Freidlin M. I., Wentzell A. D., {\it On small
%perturbations of dynamical systems}, Russian Math. Surveys 25
%(1970), 1-55.

%\bibitem{Day} Day M., {\it Mathematical Approach to the problem of noise induced exit}, Stochastic Analysis, Control, Optimization, and Applications.
%A volume in honor of W. H. Fleming. Editors: W. McEneauey, G. Yin, Q. Zhang. Birkhauser, Boston, 1999, pp 269-287.

\bibitem{Dyn} Dynkin E. B., {\it Markov Processes}, Springer-Velag, Berlin, Heidelberg, New York, 1965.

\bibitem{ET} Ekeland I., Temam R.,  {\it Convex Analysis and Variational Problems}, North Holland, Amsterdam, (1976).

\bibitem{EK86} Ethier S. N.,  Kurtz T. G, {\it  Markov processes: characterization and convergence},  Wiley Series
in Probability and Mathematical Statistics: Probability and Mathematical Statistics. John Wiley and Sons,
Inc., New York, 1986.

%\bibitem{FF} Freidlin M. I., {\it On Stochastic Perturbations of Systems with Rough Symmetry. Hierarchy of Markov Chains}, Journal of Statistical Physics, Vol. 157, No %6, pp 1031-1045, 2014.

\bibitem{FrK1}  Freidlin M., Koralov L., Wentzell A., {\it On the behavior of diffusion processes with traps}, Ann. Probab. 45 (2017), no. 5, 3202--3222.

\bibitem{FrK2}  Freidlin M., Koralov L., Wentzell A., {\it On diffusion in media with pockets of large diffusivity}, 
to appear in Probability Theory and Related Fields. 

\bibitem{FW} Freidlin M. I., Wentzell A. D., {\it Random
Perturbations of Dynamical Systems}, Springer 2012.

%\bibitem{HGB} Holmes-Cerfon M., Gortler S. J., Brenner M. P., {\it A geometrical approach to computing free-energy
%landscapes from short-ranged potentials}, Proc. Natl. Acad. Sci., 110 (1), pp E5-E14, 2013.

\bibitem{Kt} Koralov L., Sinai Y. G.,  {\it Theory of Probability and Random Processes}, 2-nd edition, Springer, 2012. 

\bibitem{Lig} Ligget T. M., {\it Continuous time Markov processes, an introduction},  Graduate Studies in Mathematics, Vol 113, AMS.

%\bibitem{K} Krylov  N.V., {\it  Nonlinear Elliptic and Parabolic Equations of the Second Order
% (Mathematics and its Applications)}, Springer 1987.
%
\bibitem{Mandl} Mandl P., {\it Analytical Treatment of One-dimensional Markov Processes}, Springer-Verlag, 1968.

\bibitem{Nov} Berlyand L., Kolpakov A., Novikov A., {\it Introduction to the network approximation method for materials modeling}, 
Encyclopedia of Mathematics and its Applications, 148. Cambridge University Press, Cambridge, 2013.

%\bibitem{MABM} Meng G., Arkus N, Brenner M. P., Manoharan V., {\it The free-energy landscape of clusters of attractive hard spheres},
%Science, 327, pp 560-563, 2010.

\bibitem{We} Wentzell A. D., {\it On lateral conditions for multidimensional diffusion processes}, Teor. Veroyatn. i Primen., 1959, Vol. 4,
no 2, pp 172–-185.

%\bibitem{LU} Ladyzenskaya O.A., Ural'zeva N.N., {\it Linear and
%Quasi-Linear Parabolic Equations} (Russian), Nauka 1967.

%\bibitem{OV} Oliveiri E., Vares M.E. {\it Large Deviations and
%Metastability}, Cambridge University Press, 2005.

\end{thebibliography}
\end{document}